\documentclass[11pt]{article}
\usepackage{amssymb,amsmath,mathrsfs,txfonts,graphicx,tikz,color}
\allowdisplaybreaks

\let\oldsection\section
\renewcommand\section{\setcounter{equation}{0}\oldsection}

\newtheorem{theorem}{\indent Theorem}[section]
\newtheorem{lemma}{\indent Lemma}[section]
\newtheorem{proposition}{\indent Proposition}[section]

\newtheorem{remark}{\indent Remark}[section]
\newtheorem{corollary}{\indent Corollary}[section]

\begin{document}

\title{\LARGE
Optimal decay rates of the compressible Euler equations
with time-dependent damping in $\mathbb R^n$: (I) under-damping case}
\author{
Shanming Ji$^{a,c}$,
Ming Mei$^{b,c,}$\thanks{Corresponding author. Emails: jism@scut.edu.cn (S.Ji), ming.mei@mcgill.ca (M.Mei)}
\\
\\
{ \small \it $^a$School of Mathematics, South China University of Technology}
\\
{ \small \it Guangzhou, Guangdong, 510641, P.~R.~China}
\\
{ \small \it $^b$Department of Mathematics, Champlain College Saint-Lambert}
\\
{ \small \it Quebec, J4P 3P2, Canada, and}
\\
{ \small \it $^c$Department of Mathematics and Statistics, McGill University}
\\
{ \small \it Montreal, Quebec, H3A 2K6, Canada}
}
\date{}

\maketitle

\begin{abstract}
This paper is concerned with the multi-dimensional compressible Euler equations
with time-dependent damping of the form $-\frac{\mu}{(1+t)^\lambda}\rho\boldsymbol u$
in $\mathbb R^n$, where $n\ge2$, $\mu>0$, and $\lambda\in[0,1)$.
When $\lambda>0$ is bigger, the damping effect time-asymptotically gets weaker, which is called under-damping.
We show the optimal decay estimates of the solutions such that
$\|\partial_x^\alpha (\rho-1)\|_{L^2(\mathbb R^n)}\approx
(1+t)^{-\frac{1+\lambda}{2}(\frac{n}{2}+|\alpha|)}$,
and $\|\partial_x^\alpha \boldsymbol u\|_{L^2(\mathbb R^n)}\approx
(1+t)^{-\frac{1+\lambda}{2}(\frac{n}{2}+|\alpha|)-\frac{1-\lambda}{2}}$,
and see how the under-damping effect influences the
structure of the Euler system. Different from the traditional view that the
stronger damping usually makes
the solutions decaying faster, here surprisingly we recognize that the
weaker damping with $0\le\lambda<1$
enhances the faster decay for the solutions. The adopted approach is
the technical Fourier analysis and the Green function method.
The main difficulties caused by the time-dependent damping lie in twofold:
non-commutativity of the Fourier transform of the linearized operator
precludes explicit expression of the fundamental solution;
time-dependent evolution implies that the Green matrix $G(t,s)$ is not
translation invariant, i.e., $G(t,s)\ne G(t-s,0)$.
We formulate the exact decay behavior of the Green matrices $G(t,s)$
with respect to $t$ and $s$ for both linear wave equations
and linear hyperbolic system, and finally derive the optimal decay
rates for the nonlinear Euler system.
\end{abstract}


\tableofcontents

\baselineskip=15pt

\section{Introduction}

\subsection{Modeling equations and background}

In this series of study, we consider the multi-dimensional compressible Euler equations with time-dependent damping
\begin{equation} \label{eq-Euler}
\begin{cases}
\displaystyle
\partial_t \rho+\nabla\cdot(\rho \boldsymbol u)=0, \\
\displaystyle
\partial_t (\rho \boldsymbol u)+\nabla\cdot(\rho \boldsymbol u\otimes \boldsymbol u)
+\nabla p(\rho)=-\frac{\mu}{(1+t)^\lambda}\rho \boldsymbol u, \\
\displaystyle
\rho|_{t=0}=\rho_0(x):=1+\tilde\rho_0(x), \quad \boldsymbol u|_{t=0}=\boldsymbol u_0(x),
\end{cases}
\end{equation}
where $x\in \mathbb R^n$, $n\ge2$, $\mu>0$, $\lambda\in[-1,1)$.
Here, the unknown functions $\rho(t,x)$ and $\boldsymbol u(t,x)$
represent the density and velocity of the fluid,
and the pressure $p(\rho)=\frac{1}{\gamma}\rho^\gamma$ with $\gamma>1$.
The initial data satisfy
\begin{equation}\label{new-1}
\rho_0(x) \to 1, \ \ \mbox{i.e.,} \ \
\tilde\rho_0(x) \to 0, \ \mbox{ and } \
\boldsymbol u_0(x) \to \boldsymbol 0 , \ \ \mbox{ as } |x|\to \infty.
\end{equation}
The damping effect of $-\frac{\mu}{(1+t)^\lambda}\rho \boldsymbol u$ is said to be under-damping for $\lambda>0$,  which is  time-asymptotically vanishing; and it is said to be over-damping for $\lambda<0$, which is time-asymptotically enhancing to $\infty$. In this paper, we are mainly interested in  the under-damping case with $\lambda \in [0,1)$, and leave the over-damping case with $\lambda \in [-1,0)$ in the second part \cite{Ji-Mei-2}.

The time-dependent damping phenomena were first proposed and studied by Wirth \cite{Wirth-JDE06,Wirth-JDE07,Wirth-MMAS04} for the linear damped wave equations, see also the
significant extension on the damped Klein-Gordon equations by Burq-Raugel-Schlag in \cite{Burq-Raugel-Schlag-2015, Burq-Raugel-Schlag-2018}, recently. Since then, the study on this subject
becomes one of hot spots, and intensively carried on, particularly, the research for Euler system involving time-dependent damping. The under- or over-damping effects with $\lambda>0$ or $\lambda<0$ makes the structure of the solutions to \eqref{eq-Euler} more complicated and various.

When $\mu=0$, the system \eqref{eq-Euler} is reduced to the pure Euler system
which usually does not possess the global-in-time solutions,
no matter how smooth the initial data are, and the singularity formed by shock waves
cannot be ignored \cite{C-D-S-W,Courant-F,Dafermos,Lax,Smoller}.

When $\mu>0$ and $\lambda=0$, the damping effect usually prevents the singularity formation
of shocks when the initial data are suitably smooth \cite{Sideris-Thomas-Wang},
but the damped solutions can still
blow up like shocks when the gradients of the initial data are big
\cite{Hailiang-Li,Wang-Chen}.
For $1$-D case, Hsiao and Liu \cite{Hsiao-Liu} first observed that the damped Euler system
is essentially equivalent to the nonlinear porous media equations, and showed
the convergence as $\|(v-{\bar v},u-{\bar u})(t)\|_{L^\infty}=O(t^{-1/2},t^{-1/2})$,
where $({\bar v}, {\bar u})(x/\sqrt{t})$ are the self-similar solutions to the corresponding
porous media equations, the so-called diffusion waves.
The relaxation-limit convergence in the weak sense was showed by Marcati and Milani in
\cite{Marcati-Milani}. After then, the convergence rates to the
diffusion waves were improved to $O(t^{-3/4},t^{-5/4})$ by Nishihara \cite{Nishihara} in
$L^2$-sense, and to $O(t^{-1},t^{-3/2})$ by Nishihara-Wang-Yang \cite{Nishihara-Wang-Yang}
in $L^1$-sense, respectively. Furthermore, Mei \cite{Mei} heuristically looked for the best
asymptotic profiles which are a kind of solutions for nonlinear diffusion equations with
certain selected initial data, and obtained much better convergence rates
$O(t^{-3/2}\ln t,t^{-2}\ln t)$.
For the multiple dimensional case, Sideris-Thomases-Wang \cite{Sideris-Thomas-Wang} first
showed the global existence of the solutions and the decay rates to the constant states as
$\|\partial_x^\alpha(\rho-1,\boldsymbol u)(t)\|_{L^2(\mathbb R^3)}
=O(t^{-\frac{3}{4}-\frac{|\alpha|}{2}}, t^{-\frac{3}{4}-\frac{|\alpha|+1}{2}})$
when the initial perturbations  are smooth enough in Sobolev space $H^l$,
which was then improved to
$O(t^{-\frac{3}{4}-\frac{|\alpha|}{2}-\frac{s}{2}},
t^{-\frac{3}{4}-\frac{|\alpha|+1}{2}-\frac{s}{2}})$ by Tan-Wu \cite{TanZ-JDE12}
for the initial data in the Besov space $H^l\cap {\dot B}^{-s}_{1,\infty}$
with $s\in [0,1]$, and to
$\|\partial_x^\alpha(\rho-1,\boldsymbol u)(t)\|_{H^{N-|\alpha|}}
=O(t^{-\frac{|\alpha|+s}{2}}, t^{-\frac{|\alpha|+s}{2}})$ by Tan-Wang \cite{TanZ-JDE13}
for the initial data in the Besov space ${\dot B}^{-s}_{2,\infty}\cap H^N$
with $s\in(0,3/2]$.
For the vaccum case, the existence of the entropy solutions and their convergence to
Barenbllat self-similar solutions were significantly studied by Huang-Pan-Wang
\cite{Huang-Pan-Wang}, Huang-Pan \cite{Huang-Pan}, Huang-Marcati-Pan
\cite{Huang-Marcati-Pan},
and Geng-Huang \cite{Geng-Huang}, respectively, and the free boundary case
with singularity was further studied by Luo-Zeng \cite{Luo-Zeng} recently.

When $\mu>0$ and $\lambda>0$, compared with the case of $\lambda=0$, the damping effect
$-\frac{\mu}{(1+t)^\lambda} \rho \boldsymbol u$ becomes weaker, we call it as under-damping.
This makes the feature of the compressible Euler system more complicated and fantastic.
For $1$-D case, Pan \cite{Pa1,Pa2} first proved that,
when $0<\lambda<1$ and the initial data around
the constant states are small enough in Sobolev space $H^1$,
then the solutions globally exist in time;
when $\lambda>1$ and the initial data are big, then the gradients of the solutions blow up
at finite time; when $\lambda=1$, the critical case, then the solutions still globally
exist for $\mu>2$, but blow up for $0<\mu\le 2$.
These results were then improved by Sugiyama \cite{Sugiyama1,Sugiyama2}
in $C^1$ space, and particularly, by Chen-Li-Li-Mei-Zhang \cite{Chen-Li-Li-Mei-Zhang} for the
global existence even with large initial data.
When the constant states at far fields are different, the convergence of the solutions
to the diffusion waves was investigated by Cui-Yin-Zhang-Zhu \cite{Cui-Yin-Zhang-Zhu}
and Li-Li-Mei-Zhang \cite{Li-Li-Mei-Zhang}, independently,
where the convergence rates obtained in \cite{Cui-Yin-Zhang-Zhu}
are better than in \cite{Li-Li-Mei-Zhang}.
In the critical case of $\lambda=1$ and $\mu>2$, by the variables scaling method for
finding the asymptotic profiles,  Geng-Lin-Mei \cite{Geng-Lin-Mei}
recognized that the roles of hyperbolicity and the damping effect for the Euler system
both are equivalently  important and cannot be ignored,  and further proved the convergence of the original solutions to the
asymptotic profiles which are artfully determined in the critical case, where the convergence rates are dependent on the physical quantity $\mu \ (>2)$.
For the multiple dimensional case $\mathbb R^n$ with $n=2,3$, Hou-Yin \cite{Hou-Yin} and
Hou-Witt-Yin \cite{Hou-Witt-Yin} first proved that, when $0<\lambda<1$ with $\mu>0$,
or $\lambda=1$ with $\mu>3-n$, once the initial data are smooth, compact supporting,
and zero-curl or not, then the solutions for the time-dependent damped Euler system
globally exist; while, when $\lambda>1$ with $\mu>0$, or $\lambda=1$ but $\mu\le 3-n$,
the solutions will blow up in finite time. The decay rates for high dimensional solutions
in the case $0<\lambda<1$ were proved by Pan \cite{PanXH-AA} very recently,
but these rates are not sufficient.

The main purpose of the present paper is to understand the structure of the solutions for
time-dependent damped Euler system as the damping effect getting weaker for $0<\lambda<1$,
and to derive the optimal decay rates of the solutions as
$\|\partial_x^\alpha (\rho-1)\|_{L^2(\mathbb R^n)}\approx
(1+t)^{-\frac{1+\lambda}{2}(\frac{n}{2}+|\alpha|)}$,
and $\|\partial_x^\alpha \boldsymbol u\|_{L^2(\mathbb R^n)}\approx
(1+t)^{-\frac{1+\lambda}{2}(\frac{n}{2}+|\alpha|)-\frac{1-\lambda}{2}}$,
by means of the technical Fourier analysis and
the Green function method.
We see from these optimal rates that the weaker damping with $0\le\lambda<1$
enhances the faster decay for the solutions. This is a bit surprise, and also subverts
the traditional view. In fact, as we show later, by taking Fourier transform to
the linearized system to
derive the fundamental solutions, we see that, when the damping is getting
less as $\lambda$ increases, the solutions in the high frequency part still decay slowly,
but the solutions in the low frequency part decay fast.

\subsection{Main results}

In order to obtain the optimal decay rates of the solutions for Euler system
\eqref{eq-Euler}, we need to build up the fundamental solutions for the corresponding
linearized system.

Let $v=\frac{2}{\gamma-1}(\sqrt{p'(\rho)}-1)
=\frac{2}{\gamma-1}(\rho^\frac{\gamma-1}{2}-1)$ and $\varpi=\frac{\gamma-1}{2}$.
Then $(v,\boldsymbol u)$ satisfies the following symmetric system
\begin{equation} \label{eq-vbdu}
\begin{cases}
\displaystyle
\partial_t v+\nabla\cdot\boldsymbol u
=-\boldsymbol u\cdot\nabla v-\varpi v\nabla\cdot \boldsymbol u,\\
\displaystyle
\partial_t \boldsymbol u+\nabla v+\frac{\mu}{(1+t)^\lambda}\boldsymbol u
=-(\boldsymbol u\cdot\nabla) \boldsymbol u-\varpi v\nabla v,\\
\displaystyle
v|_{t=0}=v_0(x), \quad \boldsymbol u|_{t=0}=\boldsymbol u_0(x),
\end{cases}
\end{equation}
where $v_0(x)=\frac{2}{\gamma-1}((1+\tilde\rho_0(x))^\frac{\gamma-1}{2}-1)$,
which behaves like $\tilde\rho_0(x)$ if the initial perturbation is small.

The optimal decay rate of the linearized system is essential
for the study of large time behavior of the time-dependent damped Euler equations.
The linearized system of \eqref{eq-vbdu} is
\begin{equation} \label{eq-vu-linear}
\begin{cases}
\displaystyle
\partial_t v+\nabla\cdot\boldsymbol u=0,\\
\displaystyle
\partial_t \boldsymbol u+\nabla v+\frac{\mu}{(1+t)^\lambda}\boldsymbol u=0,\\
\displaystyle
v|_{t=0}=v_0(x), \quad \boldsymbol u|_{t=0}=\boldsymbol u_0(x).
\end{cases}
\end{equation}
Let $u:=\Lambda^{-1}\nabla\cdot\boldsymbol u$ and
$\boldsymbol w:=\Lambda^{-1}\mathrm{curl}\,\boldsymbol u$
(with $(\mathrm{curl}\,\boldsymbol u)_j^k:=\partial_{x_j}u^k-\partial_{x_k}u^j$
for $\boldsymbol u=(u^1,\dots,u^n)$), see \cite{TanZ-JDE12} for example, where $\Lambda$ is the pseudo differential operator  defined by
$\Lambda^s v:=\mathscr{F}^{-1}(|\xi|^s\hat v(\xi))$ for $s\in\mathbb R$ (see the notations introduced below for details).
Then the linearized system \eqref{eq-vu-linear} is equivalent to
\begin{equation} \label{eq-vu}
\begin{cases}
\displaystyle
\partial_t v+\Lambda u=0,\\
\displaystyle
\partial_t u-\Lambda v+\frac{\mu}{(1+t)^\lambda} u=0,\\
\displaystyle
\partial_t \boldsymbol w+\frac{\mu}{(1+t)^\lambda} \boldsymbol w=0,\\
\displaystyle
v|_{t=0}=v_0(x), \quad u|_{t=0}=u_0(x), \quad \boldsymbol w|_{t=0}=\boldsymbol w_0(x),
\end{cases}
\end{equation}
where $u_0(x)=\Lambda^{-1}\nabla\cdot\boldsymbol u_0(x)$
and $\boldsymbol w_0(x)=\Lambda^{-1}\mathrm{curl}\,\boldsymbol u_0(x)$.
We note that the estimates on $(v,\boldsymbol u)$ are equivalent to
the estimates on $(v,u,\boldsymbol w)$ according to the relation
$$\boldsymbol u=-\Lambda^{-1}\nabla u-\Lambda^{-1}\nabla\cdot\boldsymbol w.$$
From the equation \eqref{eq-vu}$_3$, we can see that the vorticity $\boldsymbol w(t,x)$
of the linearized system decays to zero sub-exponentially
as
\[
\boldsymbol w(t,x)=\boldsymbol w_0(x)e^{-\frac{\mu}{1-\lambda}(1+t)^{1-\lambda}},
\]
which is faster than any algebraical decays.
So we only focus on the first two equations of \eqref{eq-vu}.
The Fourier transform $B(t,\xi)$ of the linear operator \eqref{eq-vu}  is time-dependent and
non-commutative (although it is diagonalizable), that is,
$B(t,\xi)B(s,\xi)\ne B(s,\xi)B(t,\xi)$ for general $s\ne t$ with
$$
B(t,\xi):=
\begin{pmatrix}
0 & -|\xi| \\
|\xi| & -\frac{\mu}{(1+t)^\lambda}
\end{pmatrix}.
$$
Therefore, the fundamental solution of the first two equations of \eqref{eq-vu}
cannot be represented as matrix exponential $e^{\int_0^tB(s,\xi)ds}$.

In order to formulate the optimal decay rates of the linearized system \eqref{eq-vu},
we consider the following two kinds of linear wave equations with time-dependent damping
\begin{equation} \label{eq-Pv}
\begin{cases}
\displaystyle
\partial_t^2 v-\Delta v+\frac{\mu}{(1+t)^\lambda} \partial_t v=0,
\quad x\in\mathbb R^n,
\\
\displaystyle
v|_{t=0}=v_1(x), \quad \partial_tv|_{t=0}=v_2(x),
\end{cases}
\end{equation}
and
\begin{equation} \label{eq-Pu}
\begin{cases}
\displaystyle
\partial_t^2 u-\Delta u+\partial_t\Big(\frac{\mu}{(1+t)^\lambda}u\Big)=0,
\quad x\in\mathbb R^n,
\\
\displaystyle
u|_{t=0}=u_1(x), \quad \partial_tu|_{t=0}=u_2(x),
\end{cases}
\end{equation}
which are satisfied by the solutions $v(t,x)$ and $u(t,x)$ of \eqref{eq-vu}, respectively.
The above two Cauchy problems \eqref{eq-Pv} and \eqref{eq-Pu} may seem similar at first glance,
but as we prove below, their optimal decay rates are totally different.
It should also be noted that the optimal decay rates derived from \eqref{eq-Pu}
are not the optimal decay rates of the solution $u$ in the linearized system \eqref{eq-vu}.
The reason is that the optimal decay rates of \eqref{eq-Pu}
are formulated with respect to arbitrary initial data $u_1(x)$ and $u_2(x)$,
while the solution $u$ in \eqref{eq-vu} corresponds to \eqref{eq-Pu}
with initial data $u_1(x)=u_0(x)$ and $u_2(x)=\Lambda v_0(x)-\mu u_0(x)$.
We will show that
there exist some cancellations between the evolution of initial data in this situation.

{\bf Notations.}
We denote $D_t=-i\partial_t$ and the $n$-dimensional Fourier transform $\mathscr{F}(v)$
of a function $v(x)$ is denoted by $\hat v(\xi)$ for simplicity.
We use $H^s=H^s(\mathbb R^n)$, $s\in\mathbb R$, to denote Sobolev spaces and
$L^p=L^p(\mathbb R^n)$, $1\le p\le \infty$, to denote the $L^p$ spaces.
The spatial derivatives $\partial_x^\alpha$ stands for
$\partial_{x_1}^{\alpha_1}\cdots\partial_{x_n}^{\alpha_n}$
with nonnegative multi-index $\alpha=(\alpha_1,\dots,\alpha_n)$, where
the order of $\alpha$ is denoted by $|\alpha|=\sum_{j=1}^{j=n}\alpha_j$, and
$\partial_x^{|\alpha|}$ stands for all the
spatial partial derivatives of order $|\alpha|$.
The pseudo differential operator $\Lambda$ is defined by
$\Lambda^s v:=\mathscr{F}^{-1}(|\xi|^s\hat v(\xi))$ for $s\in\mathbb R$.
We use $\dot H^s=\dot H^s(\mathbb R^n)$, $s\in \mathbb R$,
to denote homogeneous Sobolev spaces with
the norm $\|\cdot\|_{\dot H^s}$ defined by $\|v\|_{\dot H^s}:=\|\Lambda^sv\|_{L^2}$.
The norm $\|v\|_X^l$ stands for the $\|\cdot\|_X$ norm of the low frequency
part $v^l:=\mathscr{F}^{-1}(\chi(\xi)\hat v(\xi))$ of $v$,
while $\|v\|_X^h$ stands for the $\|\cdot\|_X$ norm of the high frequency
part $v^h:=\mathscr{F}^{-1}((1-\chi(\xi))\hat v(\xi))$ of $v$,
where $0\le\chi(\xi)\le1$ is a smooth cut-off function supported in $B_{2R}(0)$
and $\chi(\xi)\equiv1$ on $B_R(0)$ for a given $R>0$.

Throughout this paper, we also
denote $b(t)=\frac{\mu}{(1+t)^\lambda}$ with $\mu>0$ and $\lambda\in[0,1)$
and we let $C$ (or $C_j$ with $j=1,2,\dots$)
denote some positive universal constants
(may depend on the dimension $n$, the constants $\lambda$, $\mu$, $\gamma$,
and the index $\alpha$).
We use $f\lesssim g$ or $g\gtrsim f$ if $f\le Cg$
and denote $f\approx g$ if $f\lesssim g$ and $g\gtrsim f$.
For simplicity, we use $\|(f,g)\|_X$ to denote $\|f\|_X+\|g\|_X$
and $\int f:=\int_{\mathbb R^n}f(x)dx$.
The norm $\|\cdot\|_{L^2}$ will be simplified as $\|\cdot\|$ if without confusion.
For a matrix the norm $\|\cdot\|_{\max}$ is the maximum absolute value of all its elements.
We define the characteristic functions
$$
\chi_{[s\le \frac{t}{2}]}=
\chi_{[s\le \frac{t}{2}]}(s):=
\begin{cases}
1, \quad & s\le \frac{t}{2}, \\
0, \quad & \text{others},
\end{cases}
\qquad
\chi_{[s\ge \frac{t}{2}]}=
\chi_{[s\ge \frac{t}{2}]}(s):=
\begin{cases}
1, \quad & s\ge \frac{t}{2}, \\
0, \quad & \text{others}.
\end{cases}
$$
For simplicity, we denote time decay functions
\begin{equation} \label{eq-Gamma}
\Gamma(t,s):=\Big(1+(1+t)^{1+\lambda}-(1+s)^{1+\lambda}\Big)^{-\frac{1}{2}},
\qquad
\Theta(t,s):=\min\{\Gamma(t,s),(1+t)^{-\lambda}\}.
\end{equation}
There holds
$$
\Gamma(t,s)\cdot\chi_{[s\le \frac{t}{2}]}(s)\approx (1+t)^{-\frac{1+\lambda}{2}}
\approx
\Theta(t,s)\cdot\chi_{[s\le \frac{t}{2}]}(s),
\qquad
\Theta(t,s)\lesssim \Gamma(t,s).
$$

Here we always assume $\lambda\in[0,1)$ and show that
under-damping gives rise to faster decay estimates.
Our main results are stated as follows.
We present the $L^2$ and $L^q$ decay estimates of the nonlinear system \eqref{eq-vbdu}.

\begin{theorem}[Optimal $L^2$ decay rates of nonlinear Euler system] \label{th-nonlinear}
For $n\ge2$ and $\lambda\in[0,1)$, there exists a constant $\varepsilon_0>0$, such that
the solution $(v,\boldsymbol u)$ of the nonlinear system \eqref{eq-vbdu}
corresponding to initial data $(v_0,\boldsymbol u_0)$
with small energy $\|(v_0,\boldsymbol u_0)\|_{L^1\cap H^{[\frac{n}{2}]+3}}\le\varepsilon_0$
exists globally and satisfies
\begin{equation} \label{eq-nonlinear}
\begin{cases}
\|\partial_x^\alpha v\|\lesssim (1+t)^{-\frac{1+\lambda}{4}n-\frac{1+\lambda}{2}|\alpha|},
\quad &0\le |\alpha|\le [\frac{n}{2}]+1,\\
\|\partial_x^\alpha \boldsymbol u\|\lesssim
(1+t)^{-\frac{1+\lambda}{4}n-\frac{1+\lambda}{2}(|\alpha|+1)+\lambda},
\quad &0\le |\alpha|\le [\frac{n}{2}],\\
\|\partial_x^\alpha \boldsymbol u\|\lesssim
(1+t)^{-\frac{1+\lambda}{4}n-\frac{1+\lambda}{2}|\alpha|+\lambda},
\quad &|\alpha|=[\frac{n}{2}]+1, \\
\|(v,\boldsymbol u)\|_{H^{[\frac{n}{2}]+3}}\lesssim 1.
\end{cases}
\end{equation}
The first two decay estimates in \eqref{eq-nonlinear}
(i.e., the decay estimates on $\|\partial_x^\alpha v\|$
with $0\le|\alpha|\le [\frac{n}{2}]+1$ and
$\|\partial_x^\alpha \boldsymbol u\|$ with $0\le|\alpha|\le [\frac{n}{2}]$) are optimal
and consistent with the linearized system.
\end{theorem}

\begin{theorem}[Optimal $L^q$ decay estimates of nonlinear Euler system]
\label{th-nonlinear-Lp}
For $n\ge2$, $\lambda\in[0,1)$, $q\in[2,\infty]$ and $k\ge 3+[\gamma_{2,q}]$
with $\gamma_{2,q}:=n(1/2-1/q)$,
let $(v,\boldsymbol u)$ be the solution to the nonlinear system \eqref{eq-vbdu},
corresponding to the initial data $(v_0,\boldsymbol u_0)$
with small energy such that
$\|(v_0,\boldsymbol u_0)\|_{L^1\cap H^{[\frac{n}{2}]+k}}\le\varepsilon_0$,
where $\varepsilon_0>0$, is a small constant only depending on $n,q,k$ and
the constants $\gamma,\mu,\lambda$ in the system.
Then $(v,\boldsymbol u)\in L^\infty(0,+\infty;H^{[\frac{n}{2}]+k})$ and satisfies
\begin{equation} \label{eq-nonlinear-Lp}
\begin{cases}
\|\partial_x^\alpha v\|_{L^q}\lesssim
(1+t)^{-\frac{1+\lambda}{2}\gamma_{1,q}-\frac{1+\lambda}{2}|\alpha|},
\quad &0\le |\alpha|\le 1,\\
\|\boldsymbol u\|_{L^q}\lesssim
(1+t)^{-\frac{1+\lambda}{2}\gamma_{1,q}-\frac{1-\lambda}{2}},
\end{cases}
\end{equation}
where $\gamma_{1,q}=n(1-1/q)$.
All the decay estimates in \eqref{eq-nonlinear-Lp} are optimal.
\end{theorem}

For the time-dependent damped Euler equation \eqref{eq-Euler}, we have the
following decay estimates.

\begin{corollary} \label{th-Euler}
For $n\ge2$ and $\lambda\in[0,1)$, there exists a constant $\varepsilon_0>0$, such that
the solution $(\rho,\boldsymbol u)$ of the Euler equation \eqref{eq-Euler},
corresponding to the initial data $(\rho_0,\boldsymbol u_0)$
with small energy
$\|(\rho_0-1,\boldsymbol u_0)\|_{L^1\cap H^{[\frac{n}{2}]+3}}\le\varepsilon_0$,
exists globally and satisfies
\begin{equation} \label{eq-nonlinear-Euler}
\begin{cases}
\|\partial_x^\alpha (\rho-1)\|\lesssim
(1+t)^{-\frac{1+\lambda}{4}n-\frac{1+\lambda}{2}|\alpha|},
\quad &0\le |\alpha|\le [\frac{n}{2}]+1,\\
\|\partial_x^\alpha \boldsymbol u\|\lesssim
(1+t)^{-\frac{1+\lambda}{4}n-\frac{1+\lambda}{2}(|\alpha|+1)+\lambda},
\quad &0\le |\alpha|\le [\frac{n}{2}],\\
\|\partial_x^\alpha \boldsymbol u\|\lesssim
(1+t)^{-\frac{1+\lambda}{4}n-\frac{1+\lambda}{2}|\alpha|+\lambda},
\quad &|\alpha|=[\frac{n}{2}]+1, \\
\|(\rho-1,\boldsymbol u)\|_{H^{[\frac{n}{2}]+3}}\lesssim 1.
\end{cases}
\end{equation}
The first two decay estimates in \eqref{eq-nonlinear-Euler}
(i.e., the decay estimates on $\|\partial_x^\alpha (\rho-1)\|$
with $0\le|\alpha|\le [\frac{n}{2}]+1$ and
$\|\partial_x^\alpha \boldsymbol u\|$ with $0\le|\alpha|\le [\frac{n}{2}]$) are optimal.

For $n\ge2$, $\lambda\in[0,1)$, $q\in[2,\infty]$ and $k\ge 3+[\gamma_{2,q}]$
with $\gamma_{2,q}:=n(1/2-1/q)$,
let $(\rho,\boldsymbol u)$ be the solution to the Euler equation \eqref{eq-Euler}
corresponding to initial data $(\rho_0,\boldsymbol u_0)$
with small energy such that
$\|(\rho_0-1,\boldsymbol u_0)\|_{L^1\cap H^{[\frac{n}{2}]+k}}\le\varepsilon_0$,
where $\varepsilon_0>0$ is a small constant only depending on $n,q,k$ and
the constants $\gamma,\mu,\lambda$ in the system.
Then $(\rho-1,\boldsymbol u)\in L^\infty(0,+\infty;H^{[\frac{n}{2}]+k})$ and satisfies
\begin{equation} \label{eq-nonlinear-Lp-Euler}
\begin{cases}
\|\partial_x^\alpha (\rho-1)\|_{L^q}\lesssim
(1+t)^{-\frac{1+\lambda}{2}\gamma_{1,q}-\frac{1+\lambda}{2}|\alpha|},
\quad &0\le |\alpha|\le 1,\\
\|\boldsymbol u\|_{L^q}\lesssim
(1+t)^{-\frac{1+\lambda}{2}\gamma_{1,q}-\frac{1-\lambda}{2}},
\end{cases}
\end{equation}
where $\gamma_{1,q}=n(1-1/q)$.
All the decay estimates in \eqref{eq-nonlinear-Lp-Euler} are optimal.
\end{corollary}

To derive the optimal decay rates of the solutions for the Euler system with time-dependent damping
\eqref{eq-Euler},
 it is essential to investigate
the fundamental solutions to the linear system \eqref{eq-vu} and two kinds of wave equations \eqref{eq-Pv}
and \eqref{eq-Pu}. Here we state the optimal decays of the solutions for the linear wave equations \eqref{eq-Pv}
and \eqref{eq-Pu} and the linear hyperbolic system \eqref{eq-vu} as follows.

\begin{theorem}[Optimal decay rates of linear wave equations] \label{th-wave}
Let $v(t,x)$ and $u(t,x)$ be the solutions of
the Cauchy problems \eqref{eq-Pv} and \eqref{eq-Pu}
corresponding to the initial data $(v(s,x),\partial_tv(s,x))$
and $(u(s,x),\partial_tu(s,x))$ starting from the initial time $s$, respectively.
Then for $q\in[2,\infty]$ and $1\le p,r\le 2$
(or $\theta\in[0,\frac{n}{2})$), we have
\begin{align} \nonumber
\|\partial_x^\alpha v(t,\cdot)\|_{L^q}\lesssim &
\Gamma^{\gamma_{p,q}}(t,s)\cdot
\Theta^{|\alpha|}(t,s)
\\ \label{eq-optimal-v}
&\cdot\Big(\big\|(v(s,\cdot),(1+s)^\lambda\partial_tv(s,\cdot))\big\|_{L^{p}}^l
+\big\|(\partial_x^{|\alpha|+\omega_{r,q}}v(s,\cdot),(1+s)^\lambda
\partial_x^{|\alpha|-1+\omega_{r,q}}\partial_tv(s,\cdot))\big\|_{L^{r}}^h\Big),
\end{align}
and
\begin{align} \nonumber
\|\partial_x^\alpha u(t,\cdot)\|_{L^q}\lesssim &
\Big(\frac{1+t}{1+s}\Big)^\lambda\cdot
\Gamma^{\gamma_{p,q}}(t,s)\cdot
\Theta^{|\alpha|}(t,s)
\\ \label{eq-optimal-u}
&\cdot\Big(\big\|(u(s,\cdot),(1+s)^\lambda\partial_tu(s,\cdot))\big\|_{L^{p}}^l
+\big\|(\partial_x^{|\alpha|+\omega_{r,q}}u(s,\cdot),(1+s)^\lambda
\partial_x^{|\alpha|-1+\omega_{r,q}}\partial_tu(s,\cdot))\big\|_{L^{r}}^h\Big),
\end{align}
where $\gamma_{p,q}:=n(1/{p}-1/q)$
(or $\gamma_{p,q}$ replaced by $\beta_{\theta,q}:=\theta+\gamma_{2,q}$
and $\|\cdot\|_{L^{p}}$ norm replaced by $\|\cdot\|_{\dot H^{-\theta}}$),
and $\omega_{r,q}>\gamma_{r,q}$ for $(r,q)\ne(2,2)$ and $\omega_{2,2}=0$.

The decay estimates \eqref{eq-optimal-v} and \eqref{eq-optimal-u} are optimal
for all $t\ge s\ge0$ such that the ``$\lesssim$''
in \eqref{eq-optimal-v} and \eqref{eq-optimal-u}
can be replaced by ``$\approx$''
for some nontrivial initial data $(v(s,x),\partial_tv(s,x))$
and $(u(s,x),\partial_tu(s,x))$.

Moreover, there exists a number $T_0\ge0$ such that
the decay estimates \eqref{eq-optimal-v} and \eqref{eq-optimal-u} are
element-by-element optimal for $\frac{t}{2}\ge s\ge T_0$ in the following sense:
there exist four kinds of nontrivial
initial data $(v(s,x),0)$, $(0,\partial_tv(s,x))$,
$(u(s,x),0)$, and $(0,\partial_tu(s,x))$ starting from the time $s$ such that
the four corresponding solutions satisfy
\begin{align*}
\|\partial_x^\alpha v(t,\cdot)\|_{L^q}\approx &
\Gamma^{\gamma_{p,q}}(t,s)\cdot
\Theta^{|\alpha|}(t,s)
\cdot\Big(\big\|v(s,\cdot)\big\|_{L^{p}}^l
+\big\|\partial_x^{|\alpha|+\omega_{r,q}}v(s,\cdot)\big\|_{L^{r}}^h\Big),
\\
\|\partial_x^\alpha v(t,\cdot)\|_{L^q}\approx &
(1+s)^\lambda\cdot
\Gamma^{\gamma_{p,q}}(t,s)\cdot
\Theta^{|\alpha|}(t,s)
\cdot\Big(\big\|\partial_tv(s,\cdot)\big\|_{L^{p}}^l
+\big\|\partial_x^{|\alpha|-1+\omega_{r,q}}\partial_tv(s,\cdot)\big\|_{L^{r}}^h\Big),
\\
\|\partial_x^\alpha u(t,\cdot)\|_{L^q}\approx &
\Big(\frac{1+t}{1+s}\Big)^\lambda\cdot
\Gamma^{\gamma_{p,q}}(t,s)\cdot
\Theta^{|\alpha|}(t,s)
\cdot\Big(\big\|u(s,\cdot)\big\|_{L^{p}}^l
+\big\|\partial_x^{|\alpha|+\omega_{r,q}}u(s,\cdot)\big\|_{L^{r}}^h\Big),
\\
\|\partial_x^\alpha u(t,\cdot)\|_{L^q}\approx &
(1+t)^\lambda\cdot
\Gamma^{\gamma_{p,q}}(t,s)\cdot
\Theta^{|\alpha|}(t,s)
\cdot\Big(\big\|\partial_tu(s,\cdot)\big\|_{L^{p}}^l
+\big\|\partial_x^{|\alpha|-1+\omega_{r,q}}\partial_tu(s,\cdot)\big\|_{L^{r}}^h\Big),
\end{align*}
respectively.
\end{theorem}


\begin{corollary}
Let $v(t,x)$ and $u(t,x)$ be the solutions of
the Cauchy problems \eqref{eq-Pv} and \eqref{eq-Pu}
corresponding to the initial data $(v(0,x),\partial_tv(0,x))$
and $(u(0,x),\partial_tu(0,x))$ respectively.
Then for $q\in[2,\infty]$ and $1\le p,r\le 2$
(or $\theta\in[0,\frac{n}{2})$), we have
\begin{align*}
\|\partial_x^\alpha v(t,\cdot)\|_{L^q}\lesssim &
(1+t)^{-\frac{1+\lambda}{2}(\gamma_{p,q}+|\alpha|)}
\cdot\Big(\big\|(v(0,\cdot),\partial_tv(0,\cdot))\big\|_{L^{p}}^l
+\big\|(\partial_x^{|\alpha|+\omega_{r,q}}v(0,\cdot),
\partial_x^{|\alpha|-1+\omega_{r,q}}\partial_tv(0,\cdot))\big\|_{L^{r}}^h\Big),
\end{align*}
and
\begin{align*}
\|\partial_x^\alpha u(t,\cdot)\|_{L^q}\lesssim &
(1+t)^{-\frac{1+\lambda}{2}(\gamma_{p,q}+|\alpha|)+\lambda}
\cdot\Big(\big\|(u(0,\cdot),\partial_tu(s,\cdot))\big\|_{L^{p}}^l
+\big\|(\partial_x^{|\alpha|+\omega_{r,q}}u(0,\cdot),
\partial_x^{|\alpha|-1+\omega_{r,q}}\partial_tu(0,\cdot))\big\|_{L^{r}}^h\Big),
\end{align*}
where $\gamma_{p,q}:=n(1/{p}-1/q)$
(or $\gamma_{p,q}$ replaced by $\beta_{\theta,q}:=\theta+\gamma_{2,q}$
and $\|\cdot\|_{L^{p}}$ norm replaced by $\|\cdot\|_{\dot H^{-\theta}}$),
and $\omega_{r,q}>\gamma_{r,q}$ for $(r,q)\ne(2,2)$ and $\omega_{2,2}=0$.

The above decay estimates are optimal such that the ``$\lesssim$''
can be replaced by ``$\approx$''
for some nontrivial initial data $(v(0,x),\partial_tv(0,x))$
and $(u(0,x),\partial_tu(0,x))$.
\end{corollary}

\begin{remark}
The decay estimate \eqref{eq-optimal-v} for $s=0$
was first proved by Wirth \cite{Wirth-JDE07}
by developing a perfect diagonalization method.
For the application to nonlinear systems, we need to consider the
evolution of initial data starting from any $s\ge0$ to $t\ge s$
since the damping is time-dependent.
One of the main difficulties caused by the time-dependent damping is that
the evolution of the initial data starting from $s\ge0$ to $t\ge s$
is completely different from that starting from $0$ to $t-s$,
as can be seen from the estimates \eqref{eq-optimal-v} and \eqref{eq-optimal-u}.
As a consequence, the estimate on the decay rate of
$\int_0^t G(t,s)Q(s,x)ds$
is slower than $Q(t,x)$,
where $G(t,s)$ is a general Green function and $Q(t,x)$ is a general non-homogeneous term.
\end{remark}

\begin{remark}
It is surprising here that the two Cauchy problems \eqref{eq-Pv} and \eqref{eq-Pu}
decay with different rates.
We note that the function
$$\varphi(t,x):=\frac{1}{(1+t)^{\frac{1+\lambda}{2}n}}%
e^{-\frac{\mu(1+\lambda)|x|^2}{4(1+t)^{1+\lambda}}},$$
which satisfies $\frac{\mu}{(1+t)^\lambda} \partial_t \varphi=\Delta \varphi$,
is an asymptotic profile of \eqref{eq-Pv},
while $\psi(t,x):=\varphi(t,x)/(\frac{\mu}{(1+t)^\lambda})$,
which satisfies $\partial_t(\frac{\mu}{(1+t)^\lambda}\psi)=\Delta \psi$,
is a good asymptotic profile of \eqref{eq-Pu},
and $\psi(t,x)$ decays slower than $\varphi(t,x)$.
The functions $\varphi(t,x)$ and $\psi(t,x)$ decay at the same rates as
$v(t,x)$ and $u(t,x)$ proved in Theorem \ref{th-wave}.
\end{remark}

\begin{theorem}[Optimal decay rates of linear hyperbolic system] \label{th-linear}
Let $(v(t,x),u(t,x))$ be the solution of the linear hyperbolic system \eqref{eq-vu}
(the third equation of $\boldsymbol w(t,x)$ is neglected as it decays sub-exponentially)
corresponding to the initial data $(v(s,x),u(s,x))$
starting from the time $s$.
There exists a universal constant $T_0\ge0$ such that
for $q\in[2,\infty]$ and $1\le p,r\le 2$
(or $\theta\in[0,\frac{n}{2})$), and for $t\ge s\ge T_0$, we have
\begin{align} \nonumber
\|\partial_x^\alpha v(t,\cdot)\|_{L^q}
\lesssim &
\Gamma^{\gamma_{p,q}}(t,s)\cdot
\Theta^{|\alpha|}(t,s)
\cdot\Big(\big\|v(s,\cdot)\big\|_{L^{p}}^l
+\big\|\partial_x^{|\alpha|+\omega_{r,q}}v(s,\cdot)\big\|_{L^{r}}^h\Big)
\\ \label{eq-linear-v}
&+(1+s)^\lambda\cdot
\Gamma^{\gamma_{p,q}}(t,s)\cdot
\Theta^{|\alpha|+1}(t,s)
\cdot\Big(\big\| u(s,\cdot)\big\|_{L^{p}}^l
+\big\|
\partial_x^{|\alpha|+\omega_{r,q}} u(s,\cdot)\big\|_{L^{r}}^h\Big),
\end{align}
and
\begin{align} \nonumber
\|\partial_x^\alpha u(t,\cdot)\|_{L^q}
\lesssim &
\Big(\frac{1+t}{1+s}\Big)^\lambda\cdot
\Gamma^{\gamma_{p,q}}(t,s)\cdot
\Theta^{|\alpha|}(t,s)
\cdot\Big(\big\|u(s,\cdot)\big\|_{L^{p}}^l
+\big\|\partial_x^{|\alpha|+\omega_{r,q}}u(s,\cdot)\big\|_{L^{r}}^h\Big)
\\ \label{eq-linear-u}
&+(1+t)^\lambda\cdot
\Gamma^{\gamma_{p,q}}(t,s)\cdot
\Theta^{|\alpha|+1}(t,s)
\cdot\Big(\big\|v(s,\cdot)\big\|_{L^{p}}^l
+\big\|\partial_x^{|\alpha|+\omega_{r,q}} v(s,\cdot)\big\|_{L^{r}}^h\Big),
\end{align}
where $\gamma_{p,q}:=n(1/{p}-1/q)$
(or $\gamma_{p,q}$ replaced by $\beta_{\theta,q}:=\theta+\gamma_{2,q}$
and $\|\cdot\|_{L^{p}}$ norm replaced by $\|\cdot\|_{\dot H^{-\theta}}$),
and $\omega_{r,q}>\gamma_{r,q}$ for $(r,q)\ne(2,2)$ and $\omega_{2,2}=0$.

Furthermore, $u(t,\cdot)$ decays faster than \eqref{eq-linear-u}
provided one order higher regularity:
\begin{align} \nonumber
\|\partial_x^\alpha u(t,\cdot)\|_{L^q}
\lesssim &
(1+t)^\lambda\cdot
\Gamma^{\gamma_{p,q}}(t,s)\cdot
\Theta^{|\alpha|+1}(t,s)
\cdot
\Big(\big\|v(s,\cdot)\big\|_{L^{p}}^l
+\big\|\partial_x^{|\alpha|+1+\omega_{r,q}}v(s,\cdot)\big\|_{L^{r}}^h\Big)
\\ \label{eq-linear-u-opt}
&+(1+t)^\lambda(1+s)^\lambda\cdot
\Gamma^{\gamma_{p,q}}(t,s)\cdot
\Theta^{|\alpha|+2}(t,s)
\cdot
\Big(\big\| u(s,\cdot)\big\|_{L^{p}}^l
+\big\|\partial_x^{|\alpha|+1+\omega_{r,q}} u(s,\cdot)\big\|_{L^{r}}^h\Big).
\end{align}

Moreover, the decay estimates \eqref{eq-linear-v} is
element-by-element optimal for $\frac{t}{2}\ge s\ge T_0$ in the following sense:
there exist two kinds of nontrivial
initial data $(v(s,x),0)$ and $(0,u(s,x))$ starting from the time $s$ such that
the two corresponding solutions satisfy
\begin{align*}
&\|\partial_x^\alpha v(t,\cdot)\|_{L^q}
\approx
\Gamma^{\gamma_{p,q}}(t,s)\cdot
\Theta^{|\alpha|}(t,s)
\cdot\Big(\big\|v(s,\cdot)\big\|_{L^{p}}^l
+\big\|\partial_x^{|\alpha|+\omega_{r,q}}v(s,\cdot)\big\|_{L^{r}}^h\Big),
\end{align*}
and
\begin{align*}
\|\partial_x^\alpha v(t,\cdot)\|_{L^q}
\approx &
(1+s)^\lambda\cdot
\Gamma^{\gamma_{p,q}}(t,s)\cdot
\Theta^{|\alpha|+1}(t,s)
\cdot\Big(\big\| u(s,\cdot)\big\|_{L^{p}}^l
+\big\|
\partial_x^{|\alpha|+\omega_{r,q}} u(s,\cdot)\big\|_{L^{r}}^h\Big),
\end{align*}
respectively.

The decay estimate \eqref{eq-linear-u} and \eqref{eq-linear-u-opt} are
optimal with respect to $v(s,x)$ for all $\frac{t}{2}\ge s\ge T_0$
in the following sense:
there exists nontrivial initial data $(v(s,x),0)$ such that
\begin{align*}
\|\partial_x^\alpha u(t,\cdot)\|_{L^q}
\approx &
(1+t)^\lambda\cdot
\Gamma^{\gamma_{p,q}}(t,s)\cdot
\Theta^{|\alpha|+1}(t,s)
\cdot
\Big(\big\|v(s,\cdot)\big\|_{L^{p}}^l
+\big\|\partial_x^{|\alpha|+\omega_{r,q}}v(s,\cdot)\big\|_{L^{r}}^h\Big).
\end{align*}

The decay estimates \eqref{eq-linear-v} and \eqref{eq-linear-u-opt} are optimal
for all $t\ge s\ge0$ such that
\begin{align} \nonumber
\|\partial_x^\alpha v(t,\cdot)\|_{L^q}
\approx &
\Gamma^{\gamma_{p,q}}(t,s)\cdot
\Theta^{|\alpha|}(t,s)
\\ \nonumber
&\qquad\cdot\Big(\big\|v(s,\cdot)\big\|_{L^{p}}^l
+\big\|\partial_x^{|\alpha|+\omega_{r,q}}v(s,\cdot)\big\|_{L^{r}}^h
+\big\| u(s,\cdot)\big\|_{L^{p}}^l
+\big\|\partial_x^{|\alpha|+\omega_{r,q}} u(s,\cdot)\big\|_{L^{r}}^h
\Big),
\end{align}
and
\begin{align} \nonumber
\|\partial_x^\alpha u(t,\cdot)\|_{L^q}
\approx &
(1+t)^\lambda\cdot
\Gamma^{\gamma_{p,q}}(t,s)\cdot
\Theta^{|\alpha|+1}(t,s)
\\ \nonumber
&\qquad\cdot
\Big(\big\|v(s,\cdot)\big\|_{L^{p}}^l
+\big\|\partial_x^{|\alpha|+1+\omega_{r,q}}v(s,\cdot)\big\|_{L^{r}}^h
+\big\| u(s,\cdot)\big\|_{L^{p}}^l
+\big\|\partial_x^{|\alpha|+1+\omega_{r,q}} u(s,\cdot)\big\|_{L^{r}}^h
\Big),
\end{align}
for some nontrivial initial data $(v(s,x),u(s,x))$.
\end{theorem}

\begin{corollary}
Let $(v(t,x),u(t,x))$ be the solution of the linear hyperbolic system \eqref{eq-vu}
(the third equation of $\boldsymbol w(t,x)$ is neglected as it decays sub-exponentially)
corresponding to the initial data $(v(0,x),u(0,x))$.
Then for $q\in[2,\infty]$ and $1\le p,r\le 2$
(or $\theta\in[0,\frac{n}{2})$), we have
\begin{align*}
\|\partial_x^\alpha v(t,\cdot)\|_{L^q}
\approx &
(1+t)^{-\frac{1+\lambda}{2}(\gamma_{p,q}+|\alpha|)}
\\
&\qquad
\cdot\Big(\big\|v(0,\cdot)\big\|_{L^{p}}^l
+\big\|\partial_x^{|\alpha|+\omega_{r,q}}v(0,\cdot)\big\|_{L^{r}}^h
+\big\| u(0,\cdot)\big\|_{L^{p}}^l
+\big\|\partial_x^{|\alpha|+\omega_{r,q}} u(0,\cdot)\big\|_{L^{r}}^h
\Big),
\end{align*}
and
\begin{align*}
\|\partial_x^\alpha u(t,\cdot)\|_{L^q}
\approx &
(1+t)^{-\frac{1+\lambda}{2}(\gamma_{p,q}+|\alpha|)-\frac{1-\lambda}{2}}
\\
&\qquad
\cdot
\Big(\big\|v(0,\cdot)\big\|_{L^{p}}^l
+\big\|\partial_x^{|\alpha|+1+\omega_{r,q}}v(0,\cdot)\big\|_{L^{r}}^h
+\big\| u(0,\cdot)\big\|_{L^{p}}^l
+\big\|\partial_x^{|\alpha|+1+\omega_{r,q}} u(0,\cdot)\big\|_{L^{r}}^h
\Big),
\end{align*}
where $\gamma_{p,q}:=n(1/{p}-1/q)$
(or $\gamma_{p,q}$ replaced by $\beta_{\theta,q}:=\theta+\gamma_{2,q}$
and $\|\cdot\|_{L^{p}}$ norm replaced by $\|\cdot\|_{\dot H^{-\theta}}$),
and $\omega_{r,q}>\gamma_{r,q}$ for $(r,q)\ne(2,2)$ and $\omega_{2,2}=0$.
The above decay estimates are optimal.
\end{corollary}

\begin{remark}
The decay estimate \eqref{eq-linear-u} for $u$ in the linear hyperbolic system \eqref{eq-vu} with time-dependent damping
derived from the optimal decay
estimate \eqref{eq-optimal-u} in Theorem \ref{th-wave}
is not optimal,  since the initial data $u(0,x)=u_0(x)$ and
$\partial_t u(0,x)=\Lambda v_0(x)-\mu u_0(x)$ are not independent.
Cancelation occurs and the decay rate increases as in \eqref{eq-linear-u-opt}.
However, the estimate \eqref{eq-linear-u} is still of importance
in the decay estimates of the nonlinear system \eqref{eq-vbdu}
since the regularity required is one order lower than in the estimate \eqref{eq-linear-u-opt}.
\end{remark}

\begin{remark}

We would like also to note some new features and difficulties caused by the time dependent
damping of the linear system \eqref{eq-vu} and two kinds of wave equations \eqref{eq-Pv}
and \eqref{eq-Pu}.

\indent
(i) The general solutions of the
wave equation \eqref{eq-Pu} (satisfied by $u(t,x)$) decay optimally
slower than those solutions of \eqref{eq-Pv} (satisfied by $v(t,x)$);
while in the linear system \eqref{eq-vu},
$u(t,x)$ decays faster than $v(t,x)$.

\indent
(ii) The solutions to the linear system \eqref{eq-vu}
(and the linear wave equations \eqref{eq-Pv} and \eqref{eq-Pu})
decay faster as $\lambda\in[0,1)$ increases.
This may seem counterintuitive as weaker damping coefficients
give rise to solutions which decay faster.
We may understand it as follows: when $\lambda$ is larger,
the high frequencies decay slower as $e^{-C(1+t)^{1-\lambda}}$,
while the low frequencies decay faster as $e^{-C|\xi|^2(1+t)^{1+\lambda}}$,
and on the whole the increasing decay of the low frequencies
dominates the decay rate of the system,
which is faster as $\lambda$ increases.

\indent
(iii) For the application to nonlinear problems, the fundamental solution
of the linear hyperbolic system \eqref{eq-vu}
(and the linear wave equations \eqref{eq-Pv} and \eqref{eq-Pu})
starting from the time $s$ to $t$,
denoted by $G(t,s)$, is essentially dependent on $s$.
That is, $G(t,s)\ne G(t-s,0)$ since the decaying
damping $\frac{\mu}{(1+t)^\lambda}$ on $(s,t)$
is not comparable with the damping on $(0,t-s)$.

\indent
(iv) Two main difficulties occur when showing the optimal decay rates:
the first one is that we cannot express the fundamental solution
$\mathscr{E}(t,s,\xi)$ in the phase space as simply
$e^{\int_0^tB(\tau,\xi)d\tau}$ and approximated diagonalization scheme is applied
such that in the elliptic zone $Z_\mathrm{ell}^v$
$$
\mathscr{E}(t,s,\xi)=e^{\int_s^t(\sqrt{|m_v(\tau,\xi)|}
+\frac{\partial_t\sqrt{|m_v(\tau,\xi)|}}{2\sqrt{|m_v(\tau,\xi)|}})d\tau}
\tilde{\mathscr{E}}(t,s,\xi),
$$
where $\tilde{\mathscr{E}}(t,s,\xi):=MN_1(t,\xi)\mathcal{Q}(t,s,\xi)
N_1^{-1}(t,\xi)M^{-1}$, see Lemma \ref{le-DtA} below.
Therefore, we need not only to prove the lower bound of $e^{\int_s^t(\sqrt{|m_v(\tau,\xi)|}
+\frac{\partial_t\sqrt{|m_v(\tau,\xi)|}}{2\sqrt{|m_v(\tau,\xi)|}})d\tau}$,
but also to show that some elements of the matrix $\tilde{\mathscr{E}}(t,s,\xi)$
are not cancelled in the matrix product.
The other one is that the low frequencies are divided into
elliptic zone $Z_\mathrm{ell}^v$ and mixed zones,
where the frequencies in $Z_\mathrm{ell}^v$ decay slowest
but the region $Z_\mathrm{ell}^v$ is shrinking.
As a result, higher decay rates are needed for frequencies in mixed zones
in order to avoid the possible cancellations
between frequencies in different zones.
\end{remark}

The paper is organized as follows.
In Section 2 and Section 3, we formulate the optimal decay estimates of the
time-dependent damped wave equations and linear system separately.
The optimal $L^2$ and $L^q$ decay estimates of the nonlinear system
are proved in Section 4.


\section{Time-dependent damped wave equations}

We first focus on the optimal decay rates of the time-dependent damped wave equations
\eqref{eq-Pv} and \eqref{eq-Pu}.
Here we need to consider the wave equations
starting from any time $s\ge0$ to time $t\ge s$ for application to nonlinear problems,
since the evolution is not translation invariant due to the time-dependent damping.
This section is devoted to the proof of Theorem \ref{th-wave}.

Taking Fourier transforms to
the time-dependent damped wave equations
\eqref{eq-Pv} and \eqref{eq-Pu}, we have
\begin{equation} \label{eq-Pv-F}
\begin{cases}
\partial_t^2\hat v+|\xi|^2\hat v+b(t)\partial_t\hat v=0,\\
\hat v(0,\xi)=\hat v_1(\xi), \quad \partial_t\hat v(0,\xi)=\hat v_2(\xi),
\end{cases}
\end{equation}
and
\begin{equation} \label{eq-Pu-F}
\begin{cases}
\partial_t^2\hat u+|\xi|^2\hat u+\partial_t(b(t)\hat u)=0,\\
\hat u(0,\xi)=\hat u_1(\xi), \quad \partial_t\hat u(0,\xi)=\hat u_2(\xi),
\end{cases}
\end{equation}
where $b(t)=\frac{\mu}{(1+t)^\lambda}$ with $\mu>0$ and $\lambda\in[0,1)$.
The solutions can be represented in the form
\begin{align} \label{eq-Phiv}
\hat v(t,\xi)=\Phi_1^v(t,0,\xi)\hat v_1(\xi)+\Phi_2^v(t,0,\xi)\hat v_2(\xi),
\\ \label{eq-Phiu}
\hat u(t,\xi)=\Phi_1^u(t,0,\xi)\hat u_1(\xi)+\Phi_2^u(t,0,\xi)\hat u_2(\xi),
\end{align}
with Fourier multipliers $\Phi_j^v(t,s,\xi)$ and $\Phi_j^u(t,s,\xi)$, $j=1,2$,
which represent the evolution of initial data starting from $s\le t$.
A perfect diagonalization scheme was developed by Wirth \cite{Wirth-JDE06,Wirth-JDE07}
in order to handle the time-dependent operators
since the matrix is not commutative.

Let
\begin{align*}
\tilde v(t,\xi):=e^{\frac{1}{2}\int_0^tb(\tau)d\tau}\hat v(t,\xi),\\
\tilde u(t,\xi):=e^{\frac{1}{2}\int_0^tb(\tau)d\tau}\hat u(t,\xi).
\end{align*}
Then the equations in \eqref{eq-Pv-F} and \eqref{eq-Pu-F} are transformed into
\begin{align} \label{eq-Pv-tilde}
\partial_t^2\tilde v+\Big(|\xi|^2-\frac{1}{4}b^2(t)-\frac{1}{2}b'(t)\Big)\tilde v=0,
\\ \label{eq-Pu-tilde}
\partial_t^2\tilde u+\Big(|\xi|^2-\frac{1}{4}b^2(t)+\frac{1}{2}b'(t)\Big)\tilde u=0.
\end{align}
For simplicity, we denote
$$
m_v(t,\xi):=|\xi|^2-\frac{1}{4}b^2(t)-\frac{1}{2}b'(t),
\quad m_u(t,\xi):=|\xi|^2-\frac{1}{4}b^2(t)+\frac{1}{2}b'(t).
$$
One may think that the difference between $m_v(t,\xi)$ and $m_u(t,\xi)$
is of no importance since $|b'(t)|\approx \frac{1}{(1+t)^{1+\lambda}}$
is dominated by $b^2(t)\approx\frac{1}{(1+t)^{2\lambda}}$ as $\lambda\in[0,1)$.
However, we will prove that this difference makes
the solution $u(t,x)$ of \eqref{eq-Pu} essentially decay slower than
the solution $v(t,x)$ of \eqref{eq-Pv}.

We employ the diagonalization method developed by Wirth \cite{Wirth-JDE06,Wirth-JDE07}
and we pay more attention to the exact asymptotic behavior of different frequencies,
especially the low frequencies such that $m_v(t,\xi)<0$ or $m_u(t,\xi)<0$.
We need to analyze the phase-time space for both \eqref{eq-Pv-tilde} and \eqref{eq-Pu-tilde}.
For the sake of simplicity, we only write down the analysis and diagonalization of
the problem \eqref{eq-Pv-tilde} and then we highlight the difference between the two problems.
The phase-time space $(t,\xi)$ of the problem \eqref{eq-Pv-tilde}
is divided into the following parts:
\begin{align*}
Z_\mathrm{hyp}^v:&=\{(t,\xi);\sqrt{|m_v(t,\xi)|}\ge N_vb(t), m_v(t,\xi)\ge0\}, \\
Z_\mathrm{pd}^v:&=\{(t,\xi);\varepsilon_vb(t)\le
\sqrt{|m_v(t,\xi)|}\le N_vb(t), m_v(t,\xi)\ge0\}, \\
Z_\mathrm{red}^v:&=\{(t,\xi);\sqrt{|m_v(t,\xi)|}\le \varepsilon_vb(t)\}, \\
Z_\mathrm{ell}^v:&=\{(t,\xi);\sqrt{|m_v(t,\xi)|}\ge \varepsilon_vb(t),
m_v(t,\xi)\le0, t\ge t_\mathrm{ell}^v\},
\end{align*}
where $\varepsilon_v>0$ is chosen to be sufficiently small
such that the influence of the reduced zone $Z_\mathrm{red}^v$
on the fundamental solution is relatively small,
and $N_v>\varepsilon_v$, $t_\mathrm{ell}^v>0$.
There remains a bounded part
$\{(t,\xi);\sqrt{|m_v(t,\xi)|}\ge \varepsilon_vb(t),
m_v(t,\xi)\le0, t\in(0,t_\mathrm{ell}^v)\}$
which is of no influence.
The treatment of the zones,
$Z_\mathrm{hyp}^v$, $Z_\mathrm{pd}^v$, $Z_\mathrm{red}^v$, and $Z_\mathrm{ell}^v$
is similar to that in \cite{Wirth-JDE07},
here we present the treatment of the elliptic zone $Z_\mathrm{ell}^v$
in detail since this part will determine the decay rates of solutions.

For any fixed constant $c_0\ge \mu N_v$, we would call
\begin{align*}
\text{high~frequencies:}& ~(t,\xi)\in Z_\mathrm{hyp}^v, ~ |\xi|\ge c_0,\\
\text{low~frequencies:}& ~(t,\xi) \in Z_\mathrm{ell}^v, ~ \text{or~other~mixed~zones},
\end{align*}
where mixed zones are $Z_\mathrm{pd}^v$, $Z_\mathrm{red}^v$,
and $Z_\mathrm{hyp}^v$ with $|\xi|\le c_0$.
Note that the elliptic zone $Z_\mathrm{ell}^v$ is shrinking
and the frequencies in $Z_\mathrm{ell}^v$ decay slowest.

In the elliptic zone $Z_\mathrm{ell}^v$, we let $D_t:=-i\partial_t$
and $V:=(\sqrt{|m_v(t,\xi)|}\tilde v, D_t\tilde v)^\mathrm{T}$,
where $(\cdot)^\mathrm{T}$ is the transpose of a matrix or a vector.
Then the equation \eqref{eq-Pv-tilde} is converted into (note that $m_v(t,\xi)<0$)
\begin{equation} \label{eq-V}
D_tV=
\begin{pmatrix}
\frac{D_t\sqrt{|m_v(t,\xi)|}}{\sqrt{|m_v(t,\xi)|}} & \sqrt{|m_v(t,\xi)|}\\
-\sqrt{|m_v(t,\xi)|} &
\end{pmatrix}
V=:A(t,\xi)V.
\end{equation}
Let
$$
M=
\begin{pmatrix}
i & -i\\
1 & 1
\end{pmatrix}
, \quad M^{-1}=\frac{1}{2}
\begin{pmatrix}
-i & 1\\
i & 1
\end{pmatrix}.
$$
Then
\begin{equation} \label{eq-DtA}
\mathcal{D}_t-A(t,\xi)=M(\mathcal{D}_t-\mathcal{D}(t,\xi)-R(t,\xi))M^{-1},
\end{equation}
where
$$
\mathcal{D}_t=
\begin{pmatrix}
D_t & \\
 & D_t
\end{pmatrix}
, \quad \mathcal{D}(t,\xi)=
\begin{pmatrix}
-i\sqrt{|m_v(t,\xi)|} & \\
 & i\sqrt{|m_v(t,\xi)|}
\end{pmatrix}
, \quad R(t,\xi)={\textstyle\frac{D_t\sqrt{|m_v(t,\xi)|}}{2\sqrt{|m_v(t,\xi)|}}}
\begin{pmatrix}
1 & -1\\
-1 & 1
\end{pmatrix}.
$$
An important note here is that $\mathcal{D}_t\ne D_tI$
since $D_tF=-i\partial_tF$ is the time derivative of a scalar, or vector, or matrix $F$,
while $\mathcal{D}_tF$ for a matrix $F$ is a multiplier such that
$$
\mathcal{D}_tFG=D_t(FG)=(D_tF)G+F(D_tG)\ne (D_tF)G
$$
for general matrix or vector $G$.
For a vector $V$, there holds $\mathcal{D}_tV=D_t V$.

Now the matrices $\mathcal{D}_t$ and $\mathcal{D}(t,\xi)$ are diagonal
but $R(t,\xi)$ is not.
The bad thing is that $\|R(t,\xi)\|_{\max}\lesssim \frac{1}{1+t}$
(the norm $\|\cdot\|_{\max}$ for a matrix is the maximum absolute value of all its elements),
which is not uniformly bounded integrable with respect to time.
The key ingredient for the diagonalization
method developed by Wirth \cite{Wirth-JDE06,Wirth-JDE07}
is to proceed a step further, such that
\begin{equation} \label{eq-DtD}
(\mathcal{D}_t-\mathcal{D}(t,\xi)-R(t,\xi))N_1(t,\xi)=
N_1(t,\xi)(\mathcal{D}_t-\mathcal{D}(t,\xi)-F_0(t,\xi)-R_1(t,\xi)),
\end{equation}
with
$$
N^{(1)}(t,\xi)={\textstyle\frac{iD_t\sqrt{|m_v(t,\xi)|}}{2|m_v(t,\xi)|}}
\begin{pmatrix}
 & 1\\
-1 &
\end{pmatrix},
\qquad
F_0(t,\xi)={\textstyle\frac{D_t\sqrt{|m_v(t,\xi)|}}{2\sqrt{|m_v(t,\xi)|}}}
\begin{pmatrix}
1 & \\
 & 1
\end{pmatrix},
$$
and $N_1(t,\xi)=I+N^{(1)}(t,\xi)$ such that
$$
N^{(1)}(t,\xi)\mathcal{D}(t,\xi)-\mathcal{D}(t,\xi)N^{(1)}(t,\xi)=R(t,\xi)-F_0(t,\xi),
$$
and then
$$
R_1(t,\xi)=-(I+N^{(1)}(t,\xi))^{-1}(D_tN^{(1)}(t,\xi)-R(t,\xi)N^{(1)}(t,\xi)
+N^{(1)}(t,\xi)F_0(t,\xi)).
$$
Now one can verify that $\|R_1(t,\xi)\|_{\max}\lesssim \frac{1}{(1+t)^{2-\lambda}}$,
whose integral with respect to time over any interval $(s,t)$ is uniformly bounded.
We also note that
$\|N_1(t,\xi)-I\|_{\max}=\|N^{(1)}(t,\xi)\|_{\max}
\lesssim\frac{1}{(1+t)^{1-\lambda}}$
and $N_1(t,\xi)$ is uniformly bounded invertible if the $t_\mathrm{ell}^v$
in the definition of $Z_\mathrm{ell}^v$ is chosen large.

\begin{lemma} \label{le-DtA}
The fundamental solution $\mathscr{E}(t,s,\xi)$
of $\mathcal{D}_t-A(t,\xi)$ (i.e. the equation \eqref{eq-V})
for $(t,\xi)\in Z_\mathrm{ell}^v$ and $0\le s\le t$ is
\begin{align*}
\mathscr{E}(t,s,\xi)=&
MN_1(t,\xi)e^{\int_s^t(\sqrt{|m_v(\tau,\xi)|}
+\frac{\partial_t\sqrt{|m_v(\tau,\xi)|}}{2\sqrt{|m_v(\tau,\xi)|}})d\tau}\mathcal{Q}(t,s,\xi)
N_1^{-1}(t,\xi)M^{-1}
\\
=&e^{\int_s^t(\sqrt{|m_v(\tau,\xi)|}
+\frac{\partial_t\sqrt{|m_v(\tau,\xi)|}}{2\sqrt{|m_v(\tau,\xi)|}})d\tau}
\tilde{\mathscr{E}}(t,s,\xi),
\end{align*}
where $\tilde{\mathscr{E}}(t,s,\xi):=MN_1(t,\xi)\mathcal{Q}(t,s,\xi)
N_1^{-1}(t,\xi)M^{-1}$
and $\mathcal{Q}(t,s,\xi)$ is the solution of the following integral equation
\begin{equation} \label{eq-calQ}
\mathcal{Q}(t,s,\xi)=H(t,s,\xi)+
i\int_s^tH(t,\theta,\xi)R_1(\theta,\xi)\mathcal{Q}(\theta,s,\xi)d\theta,
\end{equation}
with
$$
H(t,s,\xi)=
\begin{pmatrix}
1 & 0\\
0 & e^{-2\int_s^t\sqrt{|m_v(\tau,\xi)|}d\tau}
\end{pmatrix}.
$$
Moreover,
$\|\mathcal{Q}(t,s,\xi)\|_{\max}$ is uniformly bounded
and $\|\mathcal{Q}(t,s,\xi)-H(t,s,\xi)\|_{\max}\lesssim \frac{1}{(1+s)^{1-\lambda}}$
for $(t,\xi)\in Z_\mathrm{ell}^v$ and $s\le t$.
\end{lemma}
{\it \bfseries Proof.}
According to the relation \eqref{eq-DtA} and \eqref{eq-DtD},
it suffices to prove that the fundamental solution of
$\mathcal{D}_t-\mathcal{D}(t,\xi)-F_0(t,\xi)-R_1(t,\xi)$
is $\tilde{\mathscr{E}}_0\mathcal{Q}(t,s,\xi)$
with
$$
\tilde{\mathscr{E}}_0:=e^{\int_s^t(\sqrt{|m_v(\tau,\xi)|}
+\frac{\partial_t\sqrt{|m_v(\tau,\xi)|}}{2\sqrt{|m_v(\tau,\xi)|}})d\tau}.
$$
That is, we need to show
$$
\partial_t(\tilde{\mathscr{E}}_0\mathcal{Q})=
iD_t(\tilde{\mathscr{E}}_0\mathcal{Q})
=(i\mathcal{D}(t,\xi)+iF_0(t,\xi)+iR_1(t,\xi))
(\tilde{\mathscr{E}}_0\mathcal{Q}).
$$
In fact,
$$
\partial_t(\tilde{\mathscr{E}}_0\mathcal{Q})=
(\partial_t\tilde{\mathscr{E}}_0)\mathcal{Q}
+\tilde{\mathscr{E}}_0\partial_t\mathcal{Q}
=(i\mathcal{D}+iF_0-\mathcal{H})\tilde{\mathscr{E}}_0\mathcal{Q}
+\tilde{\mathscr{E}}_0\partial_t\mathcal{Q},
$$
where
$$
\mathcal{H}(t,\xi)=i(\mathcal{D}(t,\xi)+F_0(t,\xi))
-\Big(\sqrt{|m_v(t,\xi)|}
+{\textstyle\frac{\partial_t\sqrt{|m_v(t,\xi)|}}{2\sqrt{|m_v(t,\xi)|}}}\Big)I=
\begin{pmatrix}
0 & 0 \\
0 & -2\sqrt{|m_v(t,\xi)|}
\end{pmatrix}.
$$
Noticing that $\tilde{\mathscr{E}}_0$ is scalar,
we see that $\mathcal{Q}$ is the solution of
$$
\partial_t\mathcal{Q}(t,s,\xi)=\mathcal{H}(t,\xi)\mathcal{Q}(t,s,\xi)
+iR_1(t,\xi)\mathcal{Q}(t,s,\xi),
\quad \mathcal{Q}(s,s,\xi)=I,
$$
which is equivalent to the integral equation \eqref{eq-calQ}.
As proved in Theorem 15 of \cite{Wirth-JDE07},
there holds the estimates
\begin{align*}
&\|\mathcal{Q}(t,s,\xi)-H(t,s,\xi)\|_{\max}
\\
&\lesssim
\sum_{j=1}^\infty\int_s^t\|R_1(t_1,\xi)\|_{\max}\int_s^{t_1}\|R_1(t_2,\xi)\|_{\max}
\cdots\int_s^{t_{j-1}}\|R_1(t_j,\xi)\|_{\max}dt_j\cdots dt_2dt_1 \\
&\lesssim \sum_{j=1}^\infty\frac{1}{j!}\Big(\int_s^t\|R_1(\tau,\xi)\|_{\max}d\tau\Big)^j
\lesssim e^{\int_s^t\|R_1(\tau,\xi)\|_{\max}d\tau}-1.
\end{align*}
The proof is completed since
$$\int_s^t\|R_1(\tau,\xi)\|_{\max}d\tau\lesssim \int_s^t \frac{1}{(1+\tau)^{2-\lambda}}d\tau
\lesssim \frac{1}{(1+s)^{1-\lambda}},$$
which tends to zero as $s\to\infty$.
$\hfill\Box$

The following asymptotic analysis
will be used to show the optimal decay rates
of the solutions $\hat v(t,\xi)$ and $\hat u(t,\xi)$
for equations \eqref{eq-Pv-F} and \eqref{eq-Pu-F}.

\begin{lemma} \label{le-decay-est}
For $(t,\xi)\in Z_\mathrm{ell}^v$, there holds (note that $b'(t)\le0$)
\begin{equation} \label{eq-est-v}
\begin{cases}
\displaystyle
\sqrt{|m_v(t,\xi)|}+
{\textstyle\frac{\partial_t\sqrt{|m_v(t,\xi)|}}{2\sqrt{|m_v(t,\xi)|}}}
-\frac{b(t)}{2}
\le -|\xi|^2\frac{1}{b(t)}+\frac{b'(t)}{b(t)}+|r_v(t,\xi)|,
\\[3mm] \displaystyle
\sqrt{|m_v(t,\xi)|}+
{\textstyle\frac{\partial_t\sqrt{|m_v(t,\xi)|}}{2\sqrt{|m_v(t,\xi)|}}}
-\frac{b(t)}{2}
\ge -|\xi|^2\frac{C}{b(t)}+\frac{b'(t)}{b(t)}-|r_v(t,\xi)|,
\end{cases}
\end{equation}
and for $(t,\xi)\in Z_\mathrm{ell}^u$
(the definition of zones in the phase-time space corresponding to $\tilde u$
is completely similar to that of $\tilde v$), there holds
\begin{equation} \label{eq-est-u}
\begin{cases}
\displaystyle
\sqrt{|m_u(t,\xi)|}+
{\textstyle\frac{\partial_t\sqrt{|m_u(t,\xi)|}}{2\sqrt{|m_u(t,\xi)|}}}
-\frac{b(t)}{2}
\le -|\xi|^2\frac{C_1}{b(t)}+|r_u(t,\xi)|,
\\[3mm] \displaystyle
\sqrt{|m_u(t,\xi)|}+
{\textstyle\frac{\partial_t\sqrt{|m_u(t,\xi)|}}{2\sqrt{|m_u(t,\xi)|}}}
-\frac{b(t)}{2}
\ge -|\xi|^2\frac{C_2}{b(t)}-|r_u(t,\xi)|,
\end{cases}
\end{equation}
where $|r_v(t,\xi)|\lesssim \frac{1}{(1+t)^{2-\lambda}}$ and
$|r_u(t,\xi)|\lesssim \frac{1}{(1+t)^{2-\lambda}}$
such that the integrals of $|r_v(t,\xi)|$ and $|r_u(t,\xi)|$
with respect to time are uniformly bounded.
\end{lemma}
{\it \bfseries Proof.}
Recall that
$$
m_v(t,\xi):=|\xi|^2-\frac{1}{4}b^2(t)-\frac{1}{2}b'(t),
\quad m_u(t,\xi):=|\xi|^2-\frac{1}{4}b^2(t)+\frac{1}{2}b'(t),
$$
and in the elliptic zone $Z_\mathrm{ell}^v$ or $Z_\mathrm{ell}^u$,
$m_v(t,\xi)<0$ and $\sqrt{|m_v(t,\xi)|}\ge \varepsilon_v b(t)$,
or $m_u(t,\xi)<0$ and $\sqrt{|m_u(t,\xi)|}\ge \varepsilon_u b(t)$, respectively.
Then we have
$|m_v(t,\xi)|=\frac{1}{4}b^2(t)+\frac{1}{2}b'(t)-|\xi|^2\ge \varepsilon_v^2 b^2(t)$,
$|m_v(t,\xi)|\le \frac{1}{4}b^2(t)$ and
\begin{align} \nonumber
&\sqrt{|m_v(t,\xi)|}+
{\textstyle\frac{\partial_t\sqrt{|m_v(t,\xi)|}}{2\sqrt{|m_v(t,\xi)|}}}
-\frac{b(t)}{2}
\\ \nonumber
&=\frac{|m_v(t,\xi)|^2-\frac{1}{4}b^2(t)}{\sqrt{|m_v(t,\xi)|}+\frac{b(t)}{2}}
+\frac{\frac{1}{2}b(t)b'(t)+\frac{1}{2}b''(t)}{4(\frac{1}{4}b^2(t)+\frac{1}{2}b'(t)-|\xi|^2)}
\\ \nonumber
&=\frac{-|\xi|^2}{\sqrt{|m_v(t,\xi)|}+\frac{b(t)}{2}}
+\frac{\frac{1}{2}b'(t)}{\sqrt{|m_v(t,\xi)|}+\frac{b(t)}{2}}
+\frac{\frac{1}{2}b(t)b'(t)}{4(\frac{1}{4}b^2(t)+\frac{1}{2}b'(t)-|\xi|^2)}
+\frac{\frac{1}{2}b''(t)}{4(\frac{1}{4}b^2(t)+\frac{1}{2}b'(t)-|\xi|^2)}
\\
&\le -|\xi|^2\frac{1}{b(t)}+\frac{b'(t)}{b(t)}+r_v(t,\xi),
\label{eq-zmv}
\end{align}
with
$|r_v(t,\xi)|=
\big|\frac{\frac{1}{2}b''(t)}{4(\frac{1}{4}b^2(t)+\frac{1}{2}b'(t)-|\xi|^2)}\big|
\lesssim \frac{b''(t)}{\varepsilon_v^2 b^2(t)}\lesssim \frac{1}{(1+t)^{2-\lambda}}$.
This shows the first inequality in \eqref{eq-est-v}.

As for $m_u(t,\xi)$, we have
$|m_u(t,\xi)|=\frac{1}{4}b^2(t)-\frac{1}{2}b'(t)-|\xi|^2\ge \varepsilon_u^2 b^2(t)$
and $|m_u(t,\xi)|=\frac{1}{4}b^2(t)-\frac{1}{2}b'(t)-|\xi|^2\le \frac{1}{2}b^2(t)$
since $|b'(t)|\lesssim \frac{1}{(1+t)^{1+\lambda}}$ is dominated by
$b^2(t)\approx \frac{1}{(1+t)^{2\lambda}}$
and the elliptic zone $Z_\mathrm{ell}^u$ is defined within $t\ge t_\mathrm{ell}^u$
which can be chosen large.
Now, we see that
\begin{align*}
&\sqrt{|m_u(t,\xi)|}+
{\textstyle\frac{\partial_t\sqrt{|m_u(t,\xi)|}}{2\sqrt{|m_u(t,\xi)|}}}
-\frac{b(t)}{2}
\\
&=\frac{|m_u(t,\xi)|-\frac{1}{4}b^2(t)}{\sqrt{|m_u(t,\xi)|}+\frac{b(t)}{2}}
+\frac{\frac{1}{2}b(t)b'(t)-\frac{1}{2}b''(t)}{4(\frac{1}{4}b^2(t)-\frac{1}{2}b'(t)-|\xi|^2)}
\\
&=\frac{-|\xi|^2}{\sqrt{|m_u(t,\xi)|}+\frac{b(t)}{2}}
+\frac{-\frac{1}{2}b'(t)}{\sqrt{|m_u(t,\xi)|}+\frac{b(t)}{2}}
+\frac{\frac{1}{2}b(t)b'(t)}{4(\frac{1}{4}b^2(t)-\frac{1}{2}b'(t)-|\xi|^2)}
+\frac{\frac{1}{2}b''(t)}{4(\frac{1}{4}b^2(t)-\frac{1}{2}b'(t)-|\xi|^2)}
\\
&=\frac{-|\xi|^2}{\sqrt{|m_u(t,\xi)|}+\frac{b(t)}{2}}+\bar r_u(t,\xi),
\end{align*}
with
$\frac{-|\xi|^2}{\sqrt{|m_u(t,\xi)|}+\frac{b(t)}{2}}\approx -|\xi|^2\frac{1}{b(t)}$
and
\begin{align*}
|\bar r_u(t,\xi)|=&
\Big|\frac{-\frac{1}{2}b'(t)}{\sqrt{|m_u(t,\xi)|}+\frac{b(t)}{2}}
+\frac{\frac{1}{2}b(t)b'(t)}{4(\frac{1}{4}b^2(t)-\frac{1}{2}b'(t)-|\xi|^2)}
+\frac{\frac{1}{2}b''(t)}{4(\frac{1}{4}b^2(t)-\frac{1}{2}b'(t)-|\xi|^2)}\Big|
\\
\lesssim&
\Big|\frac{-\frac{1}{2}b'(t)}{\sqrt{|m_u(t,\xi)|}+\frac{b(t)}{2}}
+\frac{\frac{1}{2}b'(t)}{b(t)}\Big|
+\Big|\frac{\frac{1}{2}b(t)b'(t)}{4(\frac{1}{4}b^2(t)-\frac{1}{2}b'(t)-|\xi|^2)}
-\frac{\frac{1}{2}b'(t)}{b(t)}\Big|
+\Big|\frac{b''(t)}{b^2(t)}\Big|
\\
\lesssim&
\frac{|b'(t)||\sqrt{|m_u(t,\xi)|}-\frac{b(t)}{2}|}%
{|b(t)||\sqrt{|m_u(t,\xi)|}+\frac{b(t)}{2}|}
+\frac{|b'(t)||2b'(t)+4|\xi|^2|}{|b(t)||m_u(t,\xi)|}
+\frac{|b''(t)|}{b^2(t)}
\\
\lesssim&
\frac{|b'(t)||-\frac{1}{2}b'(t)-|\xi|^2|}%
{|b(t)||\sqrt{|m_u(t,\xi)|}+\frac{b(t)}{2}|^2}
+\frac{|b'(t)||2b'(t)+4|\xi|^2|}{|b(t)||m_u(t,\xi)|}
+\frac{|b''(t)|}{b^2(t)}
\\
\lesssim&
|\xi|^2\frac{1}{b(t)}\cdot\frac{|b'(t)|}{b^2(t)}
+\frac{|b'(t)|^2}{b^3(t)}+\frac{|b''(t)|}{b^2(t)}.
\end{align*}
By noticing that $\frac{|b'(t)|}{b^2(t)}\lesssim \frac{1}{(1+t)^{1-\lambda}}$,
which tends to zero as $t\to\infty$,
we find that $\bar r_u(t,\xi)$ can be split into
$$\bar r_u(t,\xi)=|\xi|^2\frac{1}{b(t)}\cdot \omega(t,\xi)+r_u(t,\xi),$$
with
$$|r_u(t,\xi)|\lesssim \frac{|b'(t)|^2}{b^3(t)}+\frac{|b''(t)|}{b^2(t)}
\lesssim\frac{1}{(1+t)^{2-\lambda}},$$
and
$$\frac{-|\xi|^2}{\sqrt{|m_u(t,\xi)|}+\frac{b(t)}{2}}
+|\xi|^2\frac{1}{b(t)}\cdot \omega(t,\xi)
\approx -|\xi|^2\frac{1}{b(t)}$$
since $|\omega(t,\xi)|\lesssim \frac{1}{(1+t)^{1-\lambda}}$
and we can choose $t_\mathrm{ell}^u$ large enough
(it suffices to let $|\omega(t,\xi)|\le 1/4$).

We show that the second inequality in \eqref{eq-est-v} holds.
Note that in the $Z_\mathrm{ell}^v$,
$$\varepsilon_v^2b^2(t)\le |m_v(t,\xi)|=\frac{1}{4}b^2(t)+\frac{1}{2}b'(t)-|\xi|^2
\le \frac{1}{4}b^2(t).$$
Then \eqref{eq-zmv} reads as
\begin{align*}
&\sqrt{|m_v(t,\xi)|}+
{\textstyle\frac{\partial_t\sqrt{|m_v(t,\xi)|}}{2\sqrt{|m_v(t,\xi)|}}}
-\frac{b(t)}{2}
\\
&=\frac{-|\xi|^2}{\sqrt{|m_v(t,\xi)|}+\frac{b(t)}{2}}
+\frac{\frac{1}{2}b'(t)}{\sqrt{|m_v(t,\xi)|}+\frac{b(t)}{2}}
+\frac{\frac{1}{2}b(t)b'(t)}{4(\frac{1}{4}b^2(t)+\frac{1}{2}b'(t)-|\xi|^2)}
+\frac{\frac{1}{2}b''(t)}{4(\frac{1}{4}b^2(t)+\frac{1}{2}b'(t)-|\xi|^2)}
\\
&\ge -|\xi|^2\frac{2}{b(t)}
+\frac{\frac{1}{2}b'(t)}{\frac{b(t)}{2}+\frac{b(t)}{2}}
+\frac{\frac{1}{2}b(t)b'(t)}{b^2(t)}
+\bar r_v(t,\xi)
\\
&\ge -|\xi|^2\frac{2}{b(t)}
+\frac{b'(t)}{b(t)}
+\bar r_v(t,\xi),
\end{align*}
with
\begin{align*}
\bar r_v(t,\xi)
\ge& \Big(
\frac{\frac{1}{2}b'(t)}{\sqrt{|m_v(t,\xi)|}+\frac{b(t)}{2}}
-\frac{\frac{1}{2}b'(t)}{\frac{b(t)}{2}+\frac{b(t)}{2}}
\Big)
\\
&+\Big(
\frac{\frac{1}{2}b(t)b'(t)}{4(\frac{1}{4}b^2(t)+\frac{1}{2}b'(t)-|\xi|^2)}
-\frac{\frac{1}{2}b(t)b'(t)}{b^2(t)}
\Big)
-\Big|
\frac{\frac{1}{2}b''(t)}{4(\frac{1}{4}b^2(t)+\frac{1}{2}b'(t)-|\xi|^2)}
\Big|
\\
\gtrsim&
\frac{\frac{1}{2}b'(t)(\frac{1}{2}b(t)-\sqrt{|m_v(t,\xi)|})}%
{(\sqrt{|m_v(t,\xi)|}+\frac{b(t)}{2})b(t)}
+\frac{\frac{1}{2}b(t)b'(t)(-2b'(t)+4|\xi|^2)}%
{4(\frac{1}{4}b^2(t)+\frac{1}{2}b'(t)-|\xi|^2)b^2(t)}
-\frac{|b''(t)|}{b^2(t)}
\\
\gtrsim&
\frac{b'(t)}{b^2(t)}
\cdot \frac{-\frac{1}{2}b'(t)+|\xi|^2}{\frac{1}{2}b(t)+\sqrt{|m_v(t,\xi)|}}
-\frac{|b'(t)|^2}{b^3(t)}
+\frac{b'(t)}{b^2(t)}\cdot\frac{|\xi|^2}{b(t)}
-\frac{|b''(t)|}{b^2(t)}
\\
\gtrsim&
\frac{b'(t)}{b^2(t)}\cdot\frac{|\xi|^2}{b(t)}
-\frac{|b'(t)|^2}{b^3(t)}
-\frac{|b''(t)|}{b^2(t)}
\\
\gtrsim&
-|\xi|^2\frac{C_3}{b(t)}
-|r_v(t,\xi)|,
\end{align*}
where
$C_3:=\max_{t\ge0}\frac{|b'(t)|}{b^2(t)}\lesssim \max_{t\ge0}\frac{1}{(1+t)^{1-\lambda}}$
is bounded, and
$$
|r_v(t,\xi)|
\lesssim
\frac{|b'(t)|^2}{b^3(t)}+\frac{|b''(t)|}{b^2(t)}
\lesssim\frac{1}{(1+t)^{2-\lambda}}.
$$
Therefore, the second inequality in \eqref{eq-est-v} holds with $C=C_3+2$.
The proof is completed.
$\hfill\Box$

According to the asymptotic analysis of the frequencies, we can
formulate the following estimates.

\begin{lemma} \label{le-Phi}
The multiplies $\Phi_j^v(t,s,\xi)$ and $\Phi_j^u(t,s,\xi)$, $j=1,2$,
in the equations \eqref{eq-Phiv} and \eqref{eq-Phiu} have the following estimates:
there exist $c_0>0$, $\varepsilon\in(0,1/2)$, $C>0$, and $T_0\ge0$
(only depending on $\mu$ and $\lambda$) such that

(i) For $(t,\xi)\in Z_\mathrm{ell}^v$ and $0\le s\le t$, there hold
\begin{equation} \label{eq-Phiv-low}
|\Phi_1^v(t,s,\xi)|\lesssim e^{-C|\xi|^2\int_s^t\frac{1}{b(\tau)}d\tau},
\quad
|\Phi_2^v(t,s,\xi)|\lesssim \frac{1}{b(s)}\cdot e^{-C|\xi|^2\int_s^t\frac{1}{b(\tau)}d\tau};
\end{equation}
for $(t,\xi)\in Z_\mathrm{hyp}^v$, $0\le s\le t$, and $|\xi|\ge c_0$, there holds
$$
|\Phi_1^v(t,s,\xi)|+|\xi||\Phi_2^v(t,s,\xi)|
\lesssim e^{-(\frac{1}{2}-\varepsilon)\int_s^tb(\tau)d\tau};
$$
and for $(t,\xi)\not\in Z_\mathrm{ell}^v$ with $0\le s\le t$ and $|\xi|\le c_0$,
there hold
\begin{align*}
|\Phi_1^v(t,s,\xi)|&\lesssim e^{-C|\xi|^2\int_s^{\max\{s,t_\xi^v\}}\frac{1}{b(\tau)}d\tau
-(\frac{1}{2}-\varepsilon)\int_{\max\{s,t_\xi^v\}}^tb(\tau)d\tau},
\\
|\Phi_2^v(t,s,\xi)|&\lesssim {\textstyle\frac{1}{b(\min\{s,t_\xi^v\})}}\cdot
e^{-C|\xi|^2\int_s^{\max\{s,t_\xi^v\}}\frac{1}{b(\tau)}d\tau
-(\frac{1}{2}-\varepsilon)\int_{\max\{s,t_\xi^v\}}^tb(\tau)d\tau},
\end{align*}
where $t_\xi^v:=\sup\{t;(t,\xi)\in Z_\mathrm{ell}^v\}$.

(ii) For $(t,\xi)\in Z_\mathrm{ell}^u$ and $0\le s\le t$, there hold
\begin{equation} \label{eq-Phiu-low}
|\Phi_1^u(t,s,\xi)|\lesssim \frac{b(s)}{b(t)}\cdot
e^{-C|\xi|^2\int_s^t\frac{1}{b(\tau)}d\tau},
\quad
|\Phi_2^u(t,s,\xi)|\lesssim \frac{1}{b(t)}\cdot
e^{-C|\xi|^2\int_s^t\frac{1}{b(\tau)}d\tau};
\end{equation}
for $(t,\xi)\in Z_\mathrm{hyp}^u$, $0\le s\le t$, and $|\xi|\ge c_0$, there holds
$$
|\Phi_1^u(t,s,\xi)|+|\xi||\Phi_2^u(t,s,\xi)|
\lesssim e^{-(\frac{1}{2}-\varepsilon)\int_s^tb(\tau)d\tau};
$$
and for $(t,\xi)\not\in Z_\mathrm{ell}^u$ with $0\le s\le t$ and $|\xi|\le c_0$,
there hold
\begin{align*}
|\Phi_1^u(t,s,\xi)|&\lesssim
{\textstyle\frac{b(\min\{s,t_\xi^u\})}{b(t_\xi^u)}}\cdot
e^{-C|\xi|^2\int_s^{\max\{s,t_\xi^u\}}\frac{1}{b(\tau)}d\tau
-(\frac{1}{2}-\varepsilon)\int_{\max\{s,t_\xi^u\}}^tb(\tau)d\tau},
\\
|\Phi_2^u(t,s,\xi)|&\lesssim
{\textstyle\frac{1}{b(t_\xi^u)}}\cdot
e^{-C|\xi|^2\int_s^{\max\{s,t_\xi^u\}}\frac{1}{b(\tau)}d\tau
-(\frac{1}{2}-\varepsilon)\int_{\max\{s,t_\xi^u\}}^tb(\tau)d\tau},
\end{align*}
where $t_\xi^u:=\sup\{t;(t,\xi)\in Z_\mathrm{ell}^u\}$.

(iii) For $(t,\xi)\in Z_\mathrm{ell}^v$ and $T_0\le s\le t$,
the estimate \eqref{eq-Phiv-low} is optimal:
\begin{equation} \label{eq-Phiv-low-bk}
|\Phi_1^v(t,s,\xi)|\gtrsim
e^{-C |\xi|^2\int_s^t\frac{1}{b(\tau)}d\tau},
\quad
|\Phi_2^v(t,s,\xi)|\gtrsim
\frac{1}{b(s)}\cdot
e^{-C |\xi|^2\int_s^t\frac{1}{b(\tau)}d\tau},
\end{equation}
with another universal constant $C>0$.

(iv) For $(t,\xi)\in Z_\mathrm{ell}^u$ and $T_0\le s\le t$,
the estimate \eqref{eq-Phiu-low} is optimal:
\begin{equation} \label{eq-phiu-low-bk}
|\Phi_1^u(t,s,\xi)|\gtrsim \frac{b(s)}{b(t)}\cdot
e^{-C|\xi|^2\int_s^t\frac{1}{b(\tau)}d\tau},
\quad
|\Phi_2^u(t,s,\xi)|\gtrsim \frac{1}{b(t)}\cdot
e^{-C|\xi|^2\int_s^t\frac{1}{b(\tau)}d\tau},
\end{equation}
with another universal constant $C>0$.
\end{lemma}
{\it \bfseries Proof.}
The estimates (i) with $s=0$ was proved by Wirth in Theorem 17 of \cite{Wirth-JDE07}.
Here we need to consider $\Phi_j^v(t,s,\xi)$ with $s\le t$
for the application to nonlinear system \eqref{eq-vbdu}.
It should be noted that $\Phi_j^v(t,s,\xi)$ behaves different from
$\Phi_j^v(t-s,0,\xi)$ since the damping is time-dependent.

We first focus on the elliptic zones $Z_\mathrm{ell}^v$ and $Z_\mathrm{ell}^u$.
Using the fundamental solution ${\mathscr{E}}(t,s,\xi)$ of
$\mathcal{D}_t-A(t,\xi)$ in Lemma \ref{le-DtA},
we can express the solution of \eqref{eq-V} as
\begin{align*} 
\begin{pmatrix}
\sqrt{|m_v(t,\xi)|}\tilde v(t,\xi) \\
D_t\tilde v(t,\xi)
\end{pmatrix}
&={\mathscr{E}}(t,s,\xi)
\begin{pmatrix}
\sqrt{|m_v(s,\xi)|}\tilde v(s,\xi) \\
D_t\tilde v(s,\xi)
\end{pmatrix}
\\
&=e^{\int_s^t(\sqrt{|m_v(\tau,\xi)|}
+\frac{\partial_t\sqrt{|m_v(\tau,\xi)|}}{2\sqrt{|m_v(\tau,\xi)|}})d\tau}
\tilde{\mathscr{E}}(t,s,\xi)
\begin{pmatrix}
\sqrt{|m_v(s,\xi)|}\tilde v(s,\xi) \\
D_t\tilde v(s,\xi)
\end{pmatrix},
\end{align*}
where $\tilde{\mathscr{E}}(t,s,\xi):=MN_1(t,\xi)\mathcal{Q}(t,s,\xi)
N_1^{-1}(t,\xi)M^{-1}$
and $\|\tilde{\mathscr{E}}(t,s,\xi)\|_{\max}$ is uniformly bounded.
According to the relation
$$\tilde v(t,\xi)=e^{\frac{1}{2}\int_0^tb(\tau)d\tau}\hat v(t,\xi),$$
we arrive at (note that $D_t=-i\partial_t$)
\begin{align*} 
\begin{pmatrix}
\sqrt{|m_v(t,\xi)|}\hat v(t,\xi) \\
D_t\hat v(t,\xi)-i\frac{b(t)}{2}\hat v(t,\xi)
\end{pmatrix}
=e^{\int_s^t(\sqrt{|m_v(\tau,\xi)|}
+\frac{\partial_t\sqrt{|m_v(\tau,\xi)|}}{2\sqrt{|m_v(\tau,\xi)|}}-\frac{b(\tau)}{2})d\tau}
\cdot
\tilde{\mathscr{E}}(t,s,\xi)
\begin{pmatrix}
\sqrt{|m_v(s,\xi)|}\hat v(s,\xi) \\
D_t\hat v(s,\xi)-i\frac{b(s)}{2}\hat v(s,\xi)
\end{pmatrix}.
\end{align*}
Therefore,
\begin{align} \nonumber
\Phi_1^v(t,s,\xi)&=
\frac{1}{\sqrt{|m_v(t,\xi)|}}e^{\int_s^t(\sqrt{|m_v(\tau,\xi)|}
+\frac{\partial_t\sqrt{|m_v(\tau,\xi)|}}{2\sqrt{|m_v(\tau,\xi)|}}-\frac{b(\tau)}{2})d\tau}
(\sqrt{|m_v(s,\xi)|}[\tilde{\mathscr{E}}(t,s,\xi)]_{11}
-i\frac{b(s)}{2}[\tilde{\mathscr{E}}(t,s,\xi)]_{12}),
\\ \label{eq-zest-v}
\Phi_2^v(t,s,\xi)&=\frac{-i}{\sqrt{|m_v(t,\xi)|}}e^{\int_s^t(\sqrt{|m_v(\tau,\xi)|}
+\frac{\partial_t\sqrt{|m_v(\tau,\xi)|}}{2\sqrt{|m_v(\tau,\xi)|}}-\frac{b(\tau)}{2})d\tau}
[\tilde{\mathscr{E}}(t,s,\xi)]_{12},
\end{align}
such that
$$
\hat v(t,\xi)=\Phi_1^v(t,s,\xi)\hat v(s,\xi)+\Phi_2^v(t,s,\xi)\partial_t\hat v(s,\xi),
$$
where $[\cdot]_{jk}$ denotes the $(j,k)$-element of a matrix.
Note that in the elliptic zone $Z_\mathrm{ell}^v$, we have
$$\varepsilon_vb(t)\le\sqrt{|m_v(t,\xi)|}\le\frac{1}{2}b(t).$$
We apply the estimate \eqref{eq-est-v} in Lemma \ref{le-decay-est} to get
$$
e^{\int_s^t(\sqrt{|m_v(\tau,\xi)|}
+\frac{\partial_t\sqrt{|m_v(\tau,\xi)|}}{2\sqrt{|m_v(\tau,\xi)|}}-\frac{b(\tau)}{2})d\tau}
\lesssim
e^{-|\xi|^2\int_s^t\frac{1}{b(\tau)}d\tau}\cdot e^{\int_s^t\frac{b'(\tau)}{b(\tau)}d\tau}
=\frac{b(t)}{b(s)}\cdot e^{-|\xi|^2\int_s^t\frac{1}{b(\tau)}d\tau},
$$
which implies \eqref{eq-Phiv-low}.

Similarly, we have
\begin{align} \nonumber
\Phi_1^u(t,s,\xi)&=
\frac{1}{\sqrt{|m_u(t,\xi)|}}e^{\int_s^t(\sqrt{|m_u(\tau,\xi)|}
+\frac{\partial_t\sqrt{|m_u(\tau,\xi)|}}{2\sqrt{|m_u(\tau,\xi)|}}-\frac{b(\tau)}{2})d\tau}
(\sqrt{|m_u(s,\xi)|}[\tilde{\mathscr{E}}(t,s,\xi)]_{11}
-i\frac{b(s)}{2}[\tilde{\mathscr{E}}(t,s,\xi)]_{12}),
\\ \label{eq-zest-u}
\Phi_2^u(t,s,\xi)&=\frac{-i}{\sqrt{|m_u(t,\xi)|}}e^{\int_s^t(\sqrt{|m_u(\tau,\xi)|}
+\frac{\partial_t\sqrt{|m_u(\tau,\xi)|}}{2\sqrt{|m_u(\tau,\xi)|}}-\frac{b(\tau)}{2})d\tau}
[\tilde{\mathscr{E}}(t,s,\xi)]_{12}.
\end{align}
Here we have slightly abused the notion $\tilde{\mathscr{E}}(t,s,\xi)$,
which should be replaced by the matrix corresponding to the problem of $\tilde u(t,\xi)$.
We apply the estimate \eqref{eq-est-u} in Lemma \ref{le-decay-est} to get
$$
e^{\int_s^t(\sqrt{|m_u(\tau,\xi)|}
+\frac{\partial_t\sqrt{|m_u(\tau,\xi)|}}{2\sqrt{|m_u(\tau,\xi)|}}-\frac{b(\tau)}{2})d\tau}
\lesssim
e^{-|\xi|^2\int_s^t\frac{1}{b(\tau)}d\tau},
$$
which completes the proof of \eqref{eq-Phiu-low}.

The treatment in the zones $Z_\mathrm{hyp}^v$, $Z_\mathrm{pd}^v$, and $Z_\mathrm{red}^v$
of the phase-time space of $\tilde v(t,\xi)$
is similar to that in \cite{Wirth-JDE07}.
We note that for $(t,\xi)\in Z_\mathrm{hyp}^v$ and $|\xi|\ge c_0$,
$\frac{1}{b(s)}$, $\frac{b(s)}{b(t)}$, and $\frac{1}{b(t)}$
are all dominated by $e^{\varepsilon\int_s^tb(\tau)d\tau}$.
For $(t,\xi)\not \in Z_\mathrm{ell}^v$ and $|\xi|\ge c_0$,
we can apply the estimate \eqref{eq-Phiv-low} to $\Phi_j^v(t_\xi^v,s,\xi)$
if $s\le t_\xi^v$.
This completes the proof of (i) and the proof of (ii) follows similarly.

We prove that the estimate of $\Phi_2^v(t,s,\xi)$ in \eqref{eq-Phiv-low}
is optimal.
According to the optimal estimate \eqref{eq-est-v} in Lemma \ref{le-decay-est},
we see that for $(t,\xi)\in Z_\mathrm{ell}^v$,
$$
e^{\int_s^t(\sqrt{|m_v(\tau,\xi)|}
+\frac{\partial_t\sqrt{|m_v(\tau,\xi)|}}{2\sqrt{|m_v(\tau,\xi)|}}-\frac{b(\tau)}{2})d\tau}
\gtrsim
e^{-|\xi|^2\int_s^t\frac{C}{b(\tau)}d\tau}\cdot
e^{\int_s^t\frac{b'(\tau)}{b(\tau)}d\tau}
=\frac{b(t)}{b(s)}\cdot
e^{-|\xi|^2\int_s^t\frac{C}{b(\tau)}d\tau}.
$$
Then \eqref{eq-zest-v} reads as
\begin{align*}
|\Phi_1^v(t,s,\xi)|&=
\frac{1}{\sqrt{|m_v(t,\xi)|}}e^{\int_s^t(\sqrt{|m_v(\tau,\xi)|}
+\frac{\partial_t\sqrt{|m_v(\tau,\xi)|}}{2\sqrt{|m_v(\tau,\xi)|}}-\frac{b(\tau)}{2})d\tau}
\big|\sqrt{|m_v(s,\xi)|}[\tilde{\mathscr{E}}(t,s,\xi)]_{11}
-i\frac{b(s)}{2}[\tilde{\mathscr{E}}(t,s,\xi)]_{12}\big|
\\
&\gtrsim
e^{-|\xi|^2\int_s^t\frac{C}{b(\tau)}d\tau}
\big|
\frac{\sqrt{|m_v(s,\xi)|}}{b(s)}[\tilde{\mathscr{E}}(t,s,\xi)]_{11}
-\frac{i}{2}[\tilde{\mathscr{E}}(t,s,\xi)]_{12}
\big|,
\end{align*}
and
\begin{align*}
|\Phi_2^v(t,s,\xi)|&=
\frac{1}{\sqrt{|m_v(t,\xi)|}}e^{\int_s^t(\sqrt{|m_v(\tau,\xi)|}
+\frac{\partial_t\sqrt{|m_v(\tau,\xi)|}}{2\sqrt{|m_v(\tau,\xi)|}}-\frac{b(\tau)}{2})d\tau}
\big|[\tilde{\mathscr{E}}(t,s,\xi)]_{12}\big|,
\\
&\gtrsim
\frac{1}{b(s)}\cdot
e^{-|\xi|^2\int_s^t\frac{C}{b(\tau)}d\tau}
\big|[\tilde{\mathscr{E}}(t,s,\xi)]_{12}\big|.
\end{align*}
It suffices to show that there is no cancellation between
the elements of the matrix product of
$\tilde{\mathscr{E}}(t,s,\xi)$ such that
$\big|
\frac{\sqrt{|m_v(s,\xi)|}}{b(s)}[\tilde{\mathscr{E}}(t,s,\xi)]_{11}
-\frac{i}{2}[\tilde{\mathscr{E}}(t,s,\xi)]_{12}
\big|\gtrsim 1$ and
$|[\tilde{\mathscr{E}}(t,s,\xi)]_{12}|\gtrsim1$.
Noticing that
$$
\tilde{\mathscr{E}}(t,s,\xi)=MN_1(t,\xi)\mathcal{Q}(t,s,\xi)
N_1^{-1}(t,\xi)M^{-1}
=M(I+N^{(1)}(t,\xi))\mathcal{Q}(t,s,\xi)(I+N^{(1)}(t,\xi))^{-1}M^{-1},
$$
where $\|N^{(1)}(t,\xi)\|_{\max}\lesssim\frac{1}{(1+t)^{1-\lambda}}$
and $\|\mathcal{Q}(t,s,\xi)-H(t,s,\xi)\|_{\max}\lesssim\frac{1}{(1+s)^{1-\lambda}}$ with
$$
H(t,s,\xi)=
\begin{pmatrix}
1 & 0\\
0 & e^{-2\int_s^t\sqrt{|m_v(\tau,\xi)|}d\tau}
\end{pmatrix}
$$
as shown in Lemma \ref{le-DtA},
we can find $T_0\ge0$ such that for any $T_0\le s\le t$ and $(t,\xi)\in Z_\mathrm{ell}^v$,
there holds
$$
\|\tilde{\mathscr{E}}(t,s,\xi)
-MH(t,s,\xi)M^{-1}\|_{\max}\le \frac{1}{16},
$$
and furthermore we have
\begin{align*}
\Big\|MH(t,s,\xi)M^{-1}-
\frac{1}{2}
\begin{pmatrix}
1 & i\\
-i & 1
\end{pmatrix}
\Big\|_{\max}
=\frac{1}{2}e^{-2\int_s^t\sqrt{|m_v(\tau,\xi)|}d\tau}
\Big\|
\begin{pmatrix}
1 & -i\\
i & 1
\end{pmatrix}
\Big\|_{\max}
\le \frac{1}{2}e^{-2\int_s^t\varepsilon_vb(\tau)d\tau}
\le \frac{1}{16},
\end{align*}
if $t>s$ such that $\int_s^t\varepsilon_vb(\tau)d\tau\ge3\ln2/2$,
which is easily achieved since $\int_s^\infty b(\tau)d\tau$ is divergent.
Therefore,
$|[\tilde{\mathscr{E}}(t,s,\xi)]_{12}-\frac{i}{2}|\le \frac{1}{8}$
and $|[\tilde{\mathscr{E}}(t,s,\xi)]_{11}-\frac{1}{2}|\le \frac{1}{8}$,
which means
\begin{align*}
\Big|\big(
\frac{\sqrt{|m_v(s,\xi)|}}{b(s)}[\tilde{\mathscr{E}}(t,s,\xi)]_{11}
-\frac{i}{2}[\tilde{\mathscr{E}}(t,s,\xi)]_{12}\big)
-\big(
\frac{\sqrt{|m_v(s,\xi)|}}{b(s)}\cdot\frac{1}{2}
-\frac{i}{2}\cdot\frac{i}{2}\big)
\Big|
\le \frac{\sqrt{|m_v(s,\xi)|}}{b(s)}\frac{1}{8}
+\frac{1}{2}\frac{1}{8}\le \frac{3}{16}.
\end{align*}
It follows that
$$
\Big|\big(
\frac{\sqrt{|m_v(s,\xi)|}}{b(s)}[\tilde{\mathscr{E}}(t,s,\xi)]_{11}
-\frac{i}{2}[\tilde{\mathscr{E}}(t,s,\xi)]_{12}\big)
\Big|
\ge
\Big|\big(
\frac{\sqrt{|m_v(s,\xi)|}}{b(s)}\cdot\frac{1}{2}
-\frac{i}{2}\cdot\frac{i}{2}\big)
\Big|-\frac{3}{16}\ge\frac{1}{16},
$$
and the proof of (iii) is completed.

We turn to prove (iv) in a similar way as (iii).
According to the optimal estimate \eqref{eq-est-u} in Lemma \ref{le-decay-est},
for $(t,\xi)\in Z_\mathrm{ell}^v$, we have
$$
e^{\int_s^t(\sqrt{|m_u(\tau,\xi)|}
+\frac{\partial_t\sqrt{|m_v(\tau,\xi)|}}{2\sqrt{|m_v(\tau,\xi)|}}-\frac{b(\tau)}{2})d\tau}
\gtrsim
e^{-|\xi|^2\int_s^t\frac{C}{b(\tau)}d\tau}.
$$
Then \eqref{eq-zest-u} reads as
\begin{align*}
|\Phi_1^u(t,s,\xi)|&=
\frac{1}{\sqrt{|m_u(t,\xi)|}}e^{\int_s^t(\sqrt{|m_u(\tau,\xi)|}
+\frac{\partial_t\sqrt{|m_u(\tau,\xi)|}}{2\sqrt{|m_u(\tau,\xi)|}}-\frac{b(\tau)}{2})d\tau}
\big|\sqrt{|m_u(s,\xi)|}[\tilde{\mathscr{E}}(t,s,\xi)]_{11}
-i\frac{b(s)}{2}[\tilde{\mathscr{E}}(t,s,\xi)]_{12}\big|
\\
&\gtrsim
\frac{b(s)}{b(t)}\cdot
e^{-|\xi|^2\int_s^t\frac{C}{b(\tau)}d\tau}
\big|
\frac{\sqrt{|m_u(s,\xi)|}}{b(s)}[\tilde{\mathscr{E}}(t,s,\xi)]_{11}
-\frac{i}{2}[\tilde{\mathscr{E}}(t,s,\xi)]_{12}
\big|,
\end{align*}
and
\begin{align*}
|\Phi_2^u(t,s,\xi)|&=
\frac{1}{\sqrt{|m_u(t,\xi)|}}e^{\int_s^t(\sqrt{|m_u(\tau,\xi)|}
+\frac{\partial_t\sqrt{|m_u(\tau,\xi)|}}{2\sqrt{|m_u(\tau,\xi)|}}-\frac{b(\tau)}{2})d\tau}
\big|[\tilde{\mathscr{E}}(t,s,\xi)]_{12}\big|,
\\
&\gtrsim
\frac{1}{b(t)}\cdot
e^{-|\xi|^2\int_s^t\frac{C}{b(\tau)}d\tau}
\big|[\tilde{\mathscr{E}}(t,s,\xi)]_{12}\big|.
\end{align*}
The proof of
$|[\tilde{\mathscr{E}}(t,s,\xi)]_{12}|\gtrsim1$
and $\big|
\frac{\sqrt{|m_u(s,\xi)|}}{b(s)}[\tilde{\mathscr{E}}(t,s,\xi)]_{11}
-\frac{i}{2}[\tilde{\mathscr{E}}(t,s,\xi)]_{12}
\big|\gtrsim1$
in the case of $\tilde u(t,\xi)$ is the same as in (iii).
$\hfill\Box$

The above frequency analysis is used to show the optimal decay estimates
of the wave equations \eqref{eq-Pv} and \eqref{eq-Pu}.
Note that the time decay functions $\Gamma(t,s)$ and $\Theta(t,s)$
are defined as in \eqref{eq-Gamma}.

{\it\bfseries Proof of Theorem \ref{th-wave}.}
The estimate \eqref{eq-optimal-v} for $s=0$ was proved by Wirth \cite{Wirth-JDE07}.
Here we focus on the influence of $s$ and show that $u(t,x)$ decays slower than $v(t,x)$.
We also prove that those estimates are optimal.
According to the frequency decay estimates Lemma \ref{le-Phi}
and the representation
$$
\hat v(t,\xi)=\Phi_1^v(t,s,\xi)\hat v(s,\xi)+\Phi_2^v(t,s,\xi)\partial_t\hat v(s,\xi),
$$
we need to calculate the integral $\||\xi|^{|\alpha|}\hat v(t,\xi)\|_{L^{q'}}$
decomposed into several zones,
where $q':=p/(p-1)$ is the conjugate of $q'$ with $1'=\infty$.
For the low frequencies in the elliptic zone $(t,\xi)\in Z_\mathrm{ell}^v$,
we consider the case $p\in(1,2)$ and $q\in[2,\infty)$
and take $|\xi|^{|\alpha|}|\Phi_1^v(t,s,\xi)||\hat v(s,\xi)|$ for example.
Let $\xi_t:=\sup\{|\xi|;(t,\xi)\in Z_\mathrm{ell}^v\}$.
We have
\begin{align*}
&\int_{|\xi|\le \xi_t}
\Big(|\xi|^{|\alpha|}|\Phi_1^v(t,s,\xi)||\hat v(s,\xi)|\Big)^{q'}d\xi
\\
&\lesssim
\int_{|\xi|\le \xi_t}
|\xi|^{|\alpha|q'}
e^{-Cq'|\xi|^2\int_s^t\frac{1}{b(\tau)}d\tau}
|\hat v(s,\xi)|^{q'}d\xi
\\
&\lesssim
\Big(\int_{|\xi|\le \xi_t}
|\hat v(s,\xi)|^{p'}d\xi\Big)^{q'/{p'}}
\Big(\int_{|\xi|\le \xi_t}
|\xi|^{|\alpha|p'q'/(p'-q')}
e^{-Cp'q'/(p'-q')\cdot|\xi|^2\int_s^t\frac{1}{b(\tau)}d\tau}
d\xi\Big)^{1-q'/{p'}}
\\
&\lesssim
(\|v(s,x)\|_{L^p}^l)^{q'}
\Big(\int_0^{\xi_t}
|\xi|^{|\alpha|p'q'/(p'-q')+n-1}
e^{-Cp'q'/(p'-q')\cdot|\xi|^2\int_s^t\frac{1}{b(\tau)}d\tau}
d|\xi|\Big)^{1-q'/{p'}}
\\
&\lesssim
\big(\|v(s,x)\|_{L^p}^l
(1+(1+t)^{1+\lambda}-(1+s)^{1+\lambda})^{-\frac{1}{2}(\gamma_{p,q}+|\alpha|)}\big)^{q'},
\end{align*}
which is
\begin{equation} \label{eq-zGamma}
\||\xi|^{|\alpha|}\Phi_1^v(t,s,\xi)\hat v(s,\xi)\|_{L^q}\lesssim
\|v(s,x)\|_{L^p}^l
(1+(1+t)^{1+\lambda}-(1+s)^{1+\lambda})^{-\frac{1}{2}(\gamma_{p,q}+|\alpha|)},
\end{equation}
where we have used the fact that for general $\beta\ge0$ and $C>0$,
\begin{align*}
\int_0^\infty
|\xi|^\beta
e^{-C|\xi|^2\int_s^t\frac{1}{b(\tau)}d\tau}
d|\xi|
=&\int_0^\infty
(|\xi|^2{\textstyle\int_s^t\frac{1}{b(\tau)}d\tau})^{\frac{\beta}{2}}
e^{-C|\xi|^2\int_s^t\frac{1}{b(\tau)}d\tau}
d(|\xi|^2{\textstyle\int_s^t\frac{1}{b(\tau)}d\tau})^{\frac{1}{2}}
\cdot \Big({\int_s^t\frac{1}{b(\tau)}d\tau}\Big)^{-\frac{\beta+1}{2}}
\\
\lesssim& \Big({\int_s^t\frac{1}{b(\tau)}d\tau}\Big)^{-\frac{\beta+1}{2}}
\lesssim (1+(1+t)^{1+\lambda}-(1+s)^{1+\lambda})^{-\frac{\beta+1}{2}},
\quad t\ge s+1.
\end{align*}
We also have
$$
|\xi|^{|\alpha|}\le |\xi_t|^{|\alpha|}\lesssim b^{|\alpha|}(t)
\lesssim (1+t)^{-\lambda|\alpha|},
\quad \forall (t,\xi)\in Z_\mathrm{ell}^v,
$$
and then
\begin{align} \nonumber
\||\xi|^{|\alpha|}\Phi_1^v(t,s,\xi)\hat v(s,\xi)\|_{L^q}
&\lesssim
(1+t)^{-\lambda|\alpha|}\|\Phi_1^v(t,s,\xi)\hat v(s,\xi)\|_{L^q}
\\ \label{eq-zGamma2}
&\lesssim
\|v(s,x)\|_{L^p}^l
(1+(1+t)^{1+\lambda}-(1+s)^{1+\lambda})^{-\frac{1}{2}\gamma_{p,q}}
\cdot(1+t)^{-\lambda|\alpha|}.
\end{align}
Combining \eqref{eq-zGamma} and \eqref{eq-zGamma2} together,
we have
$$
\||\xi|^{|\alpha|}\Phi_1^v(t,s,\xi)\hat v(s,\xi)\|_{L^q}
\lesssim
\|v(s,x)\|_{L^p}^l\cdot
\Gamma^{\gamma_{p,q}}(t,s)
\cdot\Theta^{|\alpha|}(t,s).
$$

For the high frequencies such that $(t,\xi)\in Z_\mathrm{hyp}^v$ and $|\xi|\ge c_0$,
we consider the case $(r,q)\ne(2,2)$ and we have
\begin{align*}
&\int_{|\xi|\ge c_0}
\Big(|\xi|^{|\alpha|}|\Phi_1^v(t,s,\xi)||\hat v(s,\xi)|\Big)^{q'}d\xi
\\
&\lesssim
\int_{|\xi|\ge c_0}
|\xi|^{|\alpha|q'}e^{-(\frac{1}{2}-\varepsilon)q'\int_s^tb(\tau)d\tau}
|\hat v(s,\xi)|^{q'}d\xi
\\
&\lesssim
e^{-(\frac{1}{2}-\varepsilon)q'\int_s^tb(\tau)d\tau}
\Big(\int_{|\xi|\ge c_0}
\big(
|\xi|^{|\alpha|+\omega_{r,q}}
|\hat v(s,\xi)|\big)^{r'}d\xi\Big)^{q'/{r'}}
\Big(\int_{|\xi|\ge c_0}
|\xi|^{-\kappa}d\xi\Big)^{1-q'/{r'}}
\\
&\lesssim
(e^{-(\frac{1}{2}-\varepsilon)\int_s^tb(\tau)d\tau}
\|\partial_x^{|\alpha|+\omega_{r,q}}v(s,x)\|_{L^r}^h)^{q'},
\end{align*}
since $\kappa:=\omega_{r,q}r'q'/(r'-q')>n$.
Note that the sub-exponential function
$e^{-(\frac{1}{2}-\varepsilon)\int_s^tb(\tau)d\tau}$ decays faster than
$(1+(1+t)^{1+\lambda}-(1+s)^{1+\lambda})^{-\frac{1}{2}(\gamma_{p,q}+|\alpha|)}$.

For the mixed part of low frequencies such that
$(t,\xi)\not\in Z_\mathrm{ell}$ and $|\xi|\le c_0$,
we divide the proof into two cases: (i) $t_\xi^v\ge s+t_0$ and (ii) $t_\xi^v\le s+t_0$,
where $t_0\ge1$ is a constant such that
$\int_s^{s+t_0}\frac{1}{b(\tau)}d\tau\ge1$.
Note that $\frac{1}{b(\tau)}\approx (1+\tau)^\lambda$, and
$t_0$ can be chosen independent of $s$.
For case (ii) with $s<t_\xi^v\le s+t_0$, we have $|\xi|\approx b(t_\xi^v)\approx b(\tau)$
for $\tau\in(s,t_\xi^v)$, and
\begin{align*}
e^{-C_1|\xi|^2\int_s^{\max\{s,t_\xi^v\}}\frac{1}{b(\tau)}d\tau
-C_2\int_{\max\{s,t_\xi^v\}}^tb(\tau)d\tau}
\approx e^{-C_1\int_s^{\max\{s,t_\xi^v\}}{b(\tau)}d\tau
-C_2\int_{\max\{s,t_\xi^v\}}^tb(\tau)d\tau}
\lesssim
e^{-\min\{C_1,C_2\}\int_s^t{b(\tau)}d\tau},
\end{align*}
which is also true for $t_\xi^v\le s$.
As for the case (i),
we can use the following inequality for general $\beta\ge0$
\begin{align*}
&|\xi|^\beta e^{-C_1|\xi|^2\int_s^{\max\{s,t_\xi^v\}}\frac{1}{b(\tau)}d\tau
-C_2\int_{\max\{s,t_\xi^v\}}^tb(\tau)d\tau}
\\
&=\Big(|\xi|^2
{\textstyle\int_s^{t_\xi^v}\frac{1}{b(\tau)}d\tau}\Big)^{\frac{\beta}{2}}
e^{-C_1|\xi|^2\int_s^{t_\xi^v}\frac{1}{b(\tau)}d\tau}\cdot
e^{-C_2\int_{t_\xi^v}^tb(\tau)d\tau}
\Big(
{\textstyle\int_s^{t_\xi^v}\frac{1}{b(\tau)}d\tau}\Big)^{-\frac{\beta}{2}}
\\
&\lesssim
e^{-C_2\int_{t_\xi^v}^tb(\tau)d\tau}
\Big(
{\textstyle\int_s^{t_\xi^v}\frac{1}{b(\tau)}d\tau}\Big)^{-\frac{\beta}{2}}
\\
&\lesssim
e^{-C_2\int_{t_\xi^v}^tb(\tau)d\tau}
\Big(1+
{\textstyle
\int_{t_\xi^v}^t\frac{1}{b(\tau)}d\tau\Big/
\int_s^{t_\xi^v}\frac{1}{b(\tau)}d\tau}\Big)^{\frac{\beta}{2}}
\Big(
{\textstyle\int_s^{t_\xi^v}\frac{1}{b(\tau)}d\tau
+\int_{t_\xi^v}^t\frac{1}{b(\tau)}d\tau}\Big)^{-\frac{\beta}{2}}
\\
&\lesssim
\Big(
{\int_s^t\frac{1}{b(\tau)}d\tau}\Big)^{-\frac{\beta}{2}},
\end{align*}
since $\int_s^{t_\xi^v}\frac{1}{b(\tau)}d\tau
\ge\int_s^{s+t_0}\frac{1}{b(\tau)}d\tau\ge1$.
The rest of the proof is similar to the case
$(t,\xi)\in Z_\mathrm{ell}^v$.

Now we prove that the estimate \eqref{eq-optimal-v} is optimal.
The proof of the optimal decay of the estimate \eqref{eq-optimal-u}
follows in a similar way.
Without loss of generality, we assume that $s\ge1$ and $t\ge 2s$.
We show that the $L^1$-$L^q$ estimates are sharp,
other $L^p$-$L^q$ and $\dot H^{-\theta}$-$L^q$ estimates
can be deduced similarly
or using an interpolation theorem.
Let $T_0\ge0$ be the constant in Lemma \ref{le-Phi}.
If $s\ge T_0$,
we consider the initial data at the time $s$ with
$v(s,x)=0$ and $\partial_tv(s,x)=\mathscr{F}^{-1}(\chi)$
such that $\chi(\xi)$ is
a nonnegative and smooth function, $\chi(\xi)\equiv1$ for $|\xi|\le R$ and
$\mathrm{supp}\chi\subset B_{2R}(0)$.
Replacing the upper bound estimates \eqref{eq-Phiv-low}
by the optimal lower bound estimate \eqref{eq-Phiv-low-bk} of $\Phi_2^v(t,s,\xi)$
in the estimates within $Z_\mathrm{ell}^v$
shows that the frequencies in $Z_\mathrm{ell}^v$ decay not faster than the desired rates
in \eqref{eq-optimal-v}.
Note that $\hat v(s,\xi)=0$ and then $\Phi_1^v(t,s,\xi)$ has no influence.
We only need to show that the low frequencies in the mixed zones decay faster
such that the cancellation between frequencies in different zones can not happen.
In fact, $\xi_t:=\sup\{|\xi|;(t,\xi)\in Z_\mathrm{ell}^v\}\approx(1+t)^{-\lambda}$
and $t_\xi:=\sup\{t;(t,\xi)\in Z_\mathrm{ell}^v\}\approx|\xi|^{-\frac{1}{\lambda}}$,
we can estimate for $\xi_t\le|\xi|\le c_0$ and $|\xi|$ near $\xi_t$
\begin{align*}
e^{-C_1|\xi|^2\int_s^{t_\xi}\frac{1}{b(\tau)}d\tau
-C_2\int_{t_\xi}^tb(\tau)d\tau}
\lesssim e^{-C_1|\xi|^2 t_\xi^{1+\lambda}
-C_2(t^{1-\lambda}-t_\xi^{1-\lambda})}
\lesssim e^{-C_1|\xi|^{-\frac{1-\lambda}{\lambda}}-C_2(t^{1-\lambda}
-|\xi|^{-\frac{1-\lambda}{\lambda}})}
\lesssim e^{-\min\{C_1,C_2\}t^{1-\lambda}},
\end{align*}
which decays sub-exponentially and is faster than the desired decay.
We can also take the initial data $\partial_t v(s,x)=0$
and $v(s,x)=\mathscr{F}^{-1}(\chi)$,
and then using the optimal lower bound
estimate \eqref{eq-Phiv-low-bk} of $\Phi_1^v(t,s,\xi)$.

It remains to show the optimal decays for the case $s\le T_0$.
We first choose the initial data $(v(T_0,x),\partial_tv(T_0,x))$
at the time $T_0$ such that $v(t,x)$ decays
not faster than the rate in \eqref{eq-optimal-v}.
Then we consider the backward wave equation \eqref{eq-Pv}
with the initial data $(v(T_0,x),\partial_tv(T_0,x))$ at the time $T_0$ and
backward to the time $s\le T_0$.
Note that the problem is a linear wave equation with bounded damping coefficients
on a bounded time interval $(s,T_0)\subset(0,T_0)$,
and the solution remains bounded.
The proof is completed.
$\hfill\Box$

\section{Time-dependent damped linear hyperbolic system}

We next show the optimal decay estimates of the linear hyperbolic system \eqref{eq-vu}
starting from any time $s\ge0$ to time $t\ge s$ for the application to
nonlinear Euler system \eqref{eq-vbdu}.

{\it\bfseries Proof of Theorem \ref{th-linear}.}
We first prove that
\begin{align*}
\|\partial_x^\alpha v\|_{L^q}\lesssim &
\Gamma^{\gamma_{p,q}}(t,s)\cdot
\Theta^{|\alpha|+2}(t,s)
\\
&\cdot\Big(\big\|(v(s,\cdot),(1+s)^\lambda u(s,\cdot))\big\|_{L^{p}}^l
+\big\|(\partial_x^{|\alpha|+\omega_{r,q}}v(s,\cdot),(1+s)^\lambda
\partial_x^{|\alpha|+\omega_{r,q}} u(s,\cdot))\big\|_{L^{r}}^h\Big),
\end{align*}
and
\begin{align*}
\|\partial_x^\alpha u\|_{L^q}\lesssim &
\Big(\frac{1+t}{1+s}\Big)^\lambda\cdot
\Gamma^{\gamma_{p,q}}(t,s)\cdot
\Theta^{|\alpha|+2}(t,s)
\\
&\Big(\big\|(u(s,\cdot),(1+s)^\lambda v(s,\cdot))\big\|_{L^{p}}^l
+\big\|(\partial_x^{|\alpha|+\omega_{r,q}}u(s,\cdot),(1+s)^\lambda
\partial_x^{|\alpha|+\omega_{r,q}} v(s,\cdot))\big\|_{L^{r}}^h\Big),
\end{align*}
which follow from
the estimates \eqref{eq-optimal-v} and \eqref{eq-optimal-u}
in Theorem \ref{th-wave}.
That is, we regard $v(t,x)$ as a solution of \eqref{eq-Pv}
with the initial data $v(s,x)$ and $\partial_t v(s,x)=-\Lambda u(s,x)$,
and $u(t,x)$ as a solution of \eqref{eq-Pu} with the initial data
$u(s,x)$ and $\partial_t u(s,x)=\Lambda v(s,x)-b(s)u(s,x)$.
Note that,
\begin{align*}
&\big\|(1+s)^\lambda\partial_tv(s,\cdot)\big\|_{L^{p}}^l
\lesssim \big\|(1+s)^\lambda\Lambda u(s,\cdot)\big\|_{L^{p}}^l
\lesssim \big\|(1+s)^\lambda u(s,\cdot)\big\|_{L^{p}}^l,
\\
&~~\big\|(1+s)^\lambda \partial_x^{|\alpha|-1+\omega_{r,q}}
\partial_tv(s,\cdot)\big\|_{L^{r}}^h
\lesssim \big\|(1+s)^\lambda \partial_x^{|\alpha|+\omega_{r,q}}
u(s,\cdot)\big\|_{L^{r}}^h,
\end{align*}
and
\begin{align*}
&\big\|(1+s)^\lambda\partial_tu(s,\cdot)\big\|_{L^{p}}^l
\lesssim \big\|(1+s)^\lambda(\Lambda v(s,\cdot)-b(s)u(s,\cdot))\big\|_{L^{p}}^l
\lesssim \big\|(u(s,\cdot),(1+s)^\lambda v(s,\cdot))\big\|_{L^{p}}^l,
\\
&\quad\big\|(1+s)^\lambda \partial_x^{|\alpha|-1+\omega_{r,q}}
\partial_tu(s,\cdot)\big\|_{L^{r}}^h
\lesssim \big\|(\partial_x^{|\alpha|-1+\omega_{r,q}}u(s,\cdot),(1+s)^\lambda
\partial_x^{|\alpha|+\omega_{r,q}} v(s,\cdot))\big\|_{L^{r}}^h.
\end{align*}

However, the above estimates on low frequencies are not
element-by-element optimal
(the decay rate of $v(t,x)$ in dependence on the initial data $v(s,x)$ is optimal,
but that on $u(s,x)$ is not).
According to the frequency decay estimates Lemma \ref{le-Phi}
and the representation
$$
\hat v(t,\xi)=\Phi_1^v(t,s,\xi)\hat v(s,\xi)+\Phi_2^v(t,s,\xi)\partial_t\hat v(s,\xi)
=\Phi_1^v(t,s,\xi)\hat v(s,\xi)-|\xi|\Phi_2^v(t,s,\xi)\hat u(s,\xi),
$$
we can improve the decay rate of $v(t,x)$ in dependence on the initial data $u(s,x)$
by $\Theta(t,s)$
in a similar way as in the proof of Theorem \ref{th-wave}
since the decay rate is determined by the frequencies in $Z_\mathrm{ell}^v$.
This completes the proof of \eqref{eq-linear-v} and \eqref{eq-linear-u}.

We show that $u(t,x)$ decays faster than \eqref{eq-linear-u}.
According to the equation \eqref{eq-vu}$_2$, we have
\begin{equation} \label{eq-u}
u(t,x)=e^{-\int_s^tb(\tau)d\tau}u(s,x)
+\int_s^te^{-\int_\eta^tb(\tau)d\tau}\Lambda v(\eta,x)d\eta.
\end{equation}
The sub-exponential function
$e^{-\int_s^tb(\tau)d\tau}\approx e^{-C((1+t)^{1-\lambda}-(1+s)^{1-\lambda})}$
decays faster than any desired algebraical decay and
\begin{align*}
&\Big\|\partial_x^\alpha\int_s^t
e^{-\int_\eta^tb(\tau)d\tau}\Lambda v(\eta,x)d\eta\Big\|_{L^q}
\\
&\le \int_s^t e^{-\int_\eta^tb(\tau)d\tau}
\|\partial_x^\alpha\Lambda v(\eta,x)\|_{L^q}d\eta
\\
&\lesssim
\int_s^t e^{-\int_\eta^tb(\tau)d\tau}
\Gamma^{\gamma_{p,q}}(\eta,s)
\cdot\Theta^{|\alpha|+1}(\eta,s)
d\eta
\cdot\Big(\big\|v(s,\cdot)\big\|_{L^{p}}^l
+\big\|\partial_x^{|\alpha|+1+\omega_{r,q}}v(s,\cdot)\big\|_{L^{r}}^h\Big)
\\
&\ \ \ +\int_s^t e^{-\int_\eta^tb(\tau)d\tau}
(1+s)^\lambda\cdot
\Gamma^{\gamma_{p,q}}(\eta,s)
\cdot\Theta^{|\alpha|+2}(\eta,s)
d\eta
\cdot\Big(\big\|u(s,\cdot)\big\|_{L^{p}}^l
+\big\|\partial_x^{|\alpha|+1+\omega_{r,q}} u(s,\cdot)\big\|_{L^{r}}^h\Big).
\end{align*}
Integrating by parts yields
\begin{align*}
\int_s^t e^{-\int_\eta^tb(\tau)d\tau}
\Gamma^{\gamma_{p,q}}(\eta,s)
\cdot\Theta^{|\alpha|+1}(\eta,s)
d\eta
=&\int_s^t
\Gamma^{\gamma_{p,q}}(\eta,s)
\cdot\Theta^{|\alpha|+1}(\eta,s)
\frac{1}{b(\eta)}
d \big(e^{-\int_\eta^tb(\tau)d\tau}\big)
\\
\lesssim&
(1+t)^\lambda\cdot
\Gamma^{\gamma_{p,q}}(t,s)
\cdot\Theta^{|\alpha|+1}(t,s).
\end{align*}
This ends the proof of \eqref{eq-linear-u-opt}.

Let $T_0\ge0$ be the constant in Lemma \ref{le-Phi}.
We can prove that the estimates \eqref{eq-linear-v}
on $\|\partial_x^\alpha v(t,x)\|_{L^q}$
is optimal in a similar way as \eqref{eq-optimal-v} in Theorem \ref{th-wave}.
In fact, if $s\ge T_0$ we take $v(s,x)=\mathscr{F}^{-1}(\chi)$ and $u(s,x)=0$
to show the optimal decay with respect to $v(s,x)$,
such that $v(t,x)$ is a solution of \eqref{eq-Pv}
with the initial data $v(s,x)=\mathscr{F}^{-1}(\chi)$
and $\partial_t v(s,x)=-\Lambda u(s,x)=0$,
where $T_0\ge0$ is the constant in Lemma \ref{le-Phi} and $\chi(\xi)$
is the smooth function in the proof of Theorem \ref{th-wave}.
Alternatively, we take $u(s,x)=\mathscr{F}^{-1}(\chi)$ and $v(s,x)=0$
to show the optimal decay with respect to $u(s,x)$.
For the case $s\le T_0$, we apply the same procedure as in Theorem \ref{th-wave}.

Finally we show that the decay estimate \eqref{eq-linear-u-opt} is optimal
with respect to $v(s,x)$ for all $\frac{t}{2}\ge s\ge T_0$
by taking $v(s,x)=\mathscr{F}^{-1}(\chi)$ and $u(s,x)=0$.
For $(t,\xi)\in Z_\mathrm{ell}^u$ and $s\le t$,
according to \eqref{eq-zest-u} in the proof of Lemma \ref{le-Phi},
we have
\begin{eqnarray*}
\hat u(t,\xi)
&=&\Phi_1^u(t,s,\xi)\hat u(s,\xi)+\Phi_2^u(t,s,\xi)\partial_t\hat u(s,\xi) \\
&=&\Phi_2^u(t,s,\xi)|\xi|\cdot\chi(\xi)
\\
&=&\frac{-i}{\sqrt{|m_u(t,\xi)|}}e^{\int_s^t(\sqrt{|m_u(\tau,\xi)|}
+\frac{\partial_t\sqrt{|m_u(\tau,\xi)|}}{2\sqrt{|m_u(\tau,\xi)|}}-\frac{b(\tau)}{2})d\tau}
\cdot
|\xi|[\tilde{\mathscr{E}}(t,s,\xi)]_{12}\cdot\chi(\xi),
\end{eqnarray*}
and $|[\tilde{\mathscr{E}}(t,s,\xi)]_{12}|\gtrsim 1$.
The rest of the proof is similar to the proof of the optimal decay
in Theorem \ref{th-wave}.
The proof is completed.
$\hfill\Box$

\begin{remark}
The estimate \eqref{eq-linear-u} on $\|\partial_x^\alpha u(t,x)\|_{L^q}$
derived from the optimal estimate \eqref{eq-optimal-u}
is not optimal with respect to $u(s,x)$ for the linear system.
If one take $v(s,x)=0$ and $u(s,x)=\mathscr{F}^{-1}(\chi)$,
then the initial data of the wave equation satisfied by $u(s,x)$
are $u(s,x)=\mathscr{F}^{-1}(\chi)$
and $\partial_t u(s,x)=\Lambda v(s,t)-b(s)u(s,x)=-b(s)\mathscr{F}^{-1}(\chi)$.
According to the estimates in the proof of Lemma \ref{le-Phi},
we see that if $s\ge T_0$,
\begin{align*}
|\Phi_1^u(t,s,\xi)\hat u(s,\xi)|
\approx
\frac{b(s)}{b(t)}\cdot
e^{-|\xi|^2\int_s^t\frac{C}{b(\tau)}d\tau}\chi(\xi),
\end{align*}
and
\begin{align*}
|\Phi_2^v(t,s,\xi)\partial_t\hat u(s,\xi)|
\approx
\frac{1}{b(t)}\cdot
e^{-|\xi|^2\int_s^t\frac{C}{b(\tau)}d\tau}\cdot b(s)\chi(\xi)
\approx \frac{b(s)}{b(t)}\cdot
e^{-|\xi|^2\int_s^t\frac{C}{b(\tau)}d\tau}\chi(\xi).
\end{align*}
They are decaying of the same order and cancellations happen
as we can prove a faster decay \eqref{eq-linear-u-opt}.
\end{remark}

We have formulated two kinds of decay estimates on
$\|\partial_x^\alpha u(t,\cdot)\|_{L^q}$ in Theorem \ref{th-linear}:
one is \eqref{eq-linear-u} without optimal decay rates,
the other is \eqref{eq-linear-u-opt} with optimal decay rates
but the regularity required is one order higher.
In application to the nonlinear system, we can use the optimal \eqref{eq-linear-u-opt}
for the estimates of $\|\partial_x^\alpha u(t,\cdot)\|_{L^q}$ with lower index $\alpha$
and apply \eqref{eq-linear-u} to those with higher index $\alpha$.

We improve the decay estimates \eqref{eq-linear-u} on
$\|\partial_x^\alpha u(t,\cdot)\|_{L^q}$ in Theorem \ref{th-linear}
by taking advantage of the cancellation between the initial data
$u(s,x)$ and $\partial_t u(s,x)=\Lambda v(s,x)-b(s)u(s,x)$
if we regard $u(t,x)$ as a solution of the wave equation \eqref{eq-Pu}.

\begin{proposition}[Decay rates improved by cancellation]
\label{th-linear-st-can}
Let $(v(t,x),u(t,x))$ be the solution of the linear system \eqref{eq-vu}
corresponding to the initial data $(v(s,x),u(s,x))$
starting from the time $s$.
Then for $q\in[2,\infty]$ and $1\le p,r\le 2$
(or $\theta\in[0,\frac{n}{2})$),
and for $t\ge s\ge T_0$ ($T_0\ge0$ is the constant in Lemma \ref{le-Phi}), we have
\begin{align} \nonumber
\|\partial_x^\alpha u(t,\cdot)\|_{L^q}
\lesssim &
(1+t)^\lambda\cdot
\Gamma^{\gamma_{p,q}}(t,s)
\cdot\Theta^{|\alpha|+1}(t,s)
\cdot
\Big(\big\|v(s,\cdot)\big\|_{L^{p}}^l
+\big\|\partial_x^{|\alpha|+\omega_{r,q}}v(s,\cdot)\big\|_{L^{r}}^h\Big)
\\ \nonumber
&+(1+t)^\lambda(1+s)^\lambda\cdot
\Gamma^{\gamma_{p,q}}(t,s)
\cdot\Theta^{|\alpha|+2}(t,s)
\cdot
\Big(\big\| u(s,\cdot)\big\|_{L^{p}}^l
+\big\|\partial_x^{|\alpha|+\omega_{r,q}} u(s,\cdot)\big\|_{L^{r}}^h\Big)
\\ \nonumber
&+\hat C\Big(\frac{1+t}{1+s}\Big)^\lambda\cdot
\Gamma^{\gamma_{p,q}}(t,s)
\cdot\Theta^{|\alpha|}(t,s)\cdot
\Big(\frac{1}{(1+s)^{1-\lambda}}
+(1+(1+t)^{1+\lambda}-(1+s)^{1+\lambda})^{-1}\Big)
\\ \label{eq-linear-u-st-can}
&\qquad\cdot
\Big(\big\| u(s,\cdot)\big\|_{L^{p}}^l
+\big\|\partial_x^{|\alpha|+\omega_{r,q}} u(s,\cdot)\big\|_{L^{r}}^h\Big),
\end{align}
where $\hat C\ge0$ is a constant and $\gamma_{p,q}:=n(1/{p}-1/q)$
(or $\gamma_{p,q}$ replaced by $\beta_{\theta,q}:=\theta+\gamma_{2,q}$
and $\|\cdot\|_{L^{p}}$ norm replaced by $\|\cdot\|_{\dot H^{-\theta}}$),
and $\omega_{r,q}>\gamma_{r,q}$ for $(r,q)\ne(2,2)$ and $\omega_{2,2}=0$.
The decay estimate \eqref{eq-linear-u-st-can} is optimal
with respect to $v(s,x)$ for all $\frac{t}{2}\ge s\ge T_0$.
\end{proposition}
{\it\bfseries Proof.}
If $\hat C=0$, the decay rates in \eqref{eq-linear-u-st-can} are equal to
that in \eqref{eq-linear-u-opt} in Theorem \ref{th-linear},
but the regularity required is one order lower.
We note that the estimates on $\|\partial_x^\alpha u(t,\cdot)\|_{L^q}$ in
\eqref{eq-linear-u-opt} are deduced from the
optimal decay estimates on $\|\nabla\partial_x^\alpha v(t,\cdot)\|_{L^q}$,
which requires regularity one order higher.
Noticing that cancellations happen in the evolution between the initial data
if we regard $u(t,x)$ as a solution of the wave equation \eqref{eq-Pu},
we make advantage of the cancellation to improve the decay estimates
without the one order higher regularity.

Similar to the proof of Lemma \ref{le-Phi} but with more precise estimates
concerned with the possible cancellations,
for $(t,\xi)\in Z_\mathrm{ell}^u$ and $s\le t$,
we have
\begin{eqnarray*}
\hat u(t,\xi)
&=&\Phi_1^u(t,s,\xi)\hat u(s,\xi)+\Phi_2^u(t,s,\xi)\partial_t\hat u(s,\xi)
\\
&=&\Phi_1^u(t,s,\xi)\hat u(s,\xi)+\Phi_2^u(t,s,\xi)(|\xi|\hat v(s,\xi)-b(s)\hat u(s,\xi))
\\
&=&\big(\Phi_1^u(t,s,\xi)-b(s)\Phi_2^u(t,s,\xi)\big)\hat u(s,\xi)
+|\xi|\Phi_2^u(t,s,\xi)\hat v(s,\xi)
\\
&=&\frac{1}{\sqrt{|m_u(t,\xi)|}}e^{\int_s^t(\sqrt{|m_u(\tau,\xi)|}
+\frac{\partial_t\sqrt{|m_u(\tau,\xi)|}}{2\sqrt{|m_u(\tau,\xi)|}}-\frac{b(\tau)}{2})d\tau}
\\
&&\ \ \ \cdot\Big(
(\sqrt{|m_u(s,\xi)|}[\tilde{\mathscr{E}}(t,s,\xi)]_{11}
+i\frac{b(s)}{2}[\tilde{\mathscr{E}}(t,s,\xi)]_{12})\hat u(s,\xi)
-i|\xi|[\tilde{\mathscr{E}}(t,s,\xi)]_{12}\hat v(s,\xi)
\Big),
\end{eqnarray*}
according to \eqref{eq-zest-u} in the proof of Lemma \ref{le-Phi},
where we have proved that there are no cancelations between
$$\sqrt{|m_u(s,\xi)|}[\tilde{\mathscr{E}}(t,s,\xi)]_{11}
-i\frac{b(s)}{2}[\tilde{\mathscr{E}}(t,s,\xi)]_{12},$$
and here we show that the leading terms within the summation
$$\sqrt{|m_u(s,\xi)|}[\tilde{\mathscr{E}}(t,s,\xi)]_{11}
+i\frac{b(s)}{2}[\tilde{\mathscr{E}}(t,s,\xi)]_{12}$$
cancel each other.
In fact, noticing that
$$
\tilde{\mathscr{E}}(t,s,\xi)=MN_1(t,\xi)\mathcal{Q}(t,s,\xi)
N_1^{-1}(t,\xi)M^{-1}
=M(I+N^{(1)}(t,\xi))\mathcal{Q}(t,s,\xi)(I+N^{(1)}(t,\xi))^{-1}M^{-1},
$$
where $\|N^{(1)}(t,\xi)\|_{\max}\lesssim\frac{1}{(1+t)^{1-\lambda}}$
and $\|\mathcal{Q}(t,s,\xi)-H(t,s,\xi)\|_{\max}\lesssim\frac{1}{(1+s)^{1-\lambda}}$ with
$$
H(t,s,\xi)=
\begin{pmatrix}
1 & 0\\
0 & e^{-2\int_s^t\sqrt{|m_u(\tau,\xi)|}d\tau}
\end{pmatrix}
$$
as shown in Lemma \ref{le-DtA},
we have
\begin{equation} \label{eq-zcancel}
\|\tilde{\mathscr{E}}(t,s,\xi)
-MH(t,s,\xi)M^{-1}\|_{\max}
\lesssim
\frac{1}{(1+s)^{1-\lambda}},
\end{equation}
and
\begin{equation} \label{eq-zcancel-a}
\Big\|MH(t,s,\xi)M^{-1}-
\frac{1}{2}
\begin{pmatrix}
1 & i\\
-i & 1
\end{pmatrix}
\Big\|_{\max}
=\frac{1}{2}e^{-2\int_s^t\sqrt{|m_u(\tau,\xi)|}d\tau}
\Big\|
\begin{pmatrix}
1 & -i\\
i & 1
\end{pmatrix}
\Big\|_{\max}
\le \frac{1}{2}e^{-2\int_s^t\varepsilon_ub(\tau)d\tau}.
\end{equation}
Therefore, \eqref{eq-zcancel} and \eqref{eq-zcancel-a} imply
\begin{equation} \label{eq-zapprox}
\Big\|\tilde{\mathscr{E}}(t,s,\xi)-
\frac{1}{2}
\begin{pmatrix}
1 & i\\
-i & 1
\end{pmatrix}
\Big\|_{\max}
\lesssim
\frac{1}{(1+s)^{1-\lambda}}
+e^{-2\int_s^t\varepsilon_ub(\tau)d\tau},
\end{equation}
which means
\begin{align*}
\Big|
\frac{\sqrt{|m_u(s,\xi)|}}{b(s)}[\tilde{\mathscr{E}}(t,s,\xi)]_{11}
+\frac{i}{2}[\tilde{\mathscr{E}}(t,s,\xi)]_{12}\Big|
\lesssim&
\Big|
\frac{\sqrt{|m_u(s,\xi)|}}{b(s)}\cdot\frac{1}{2}
+\frac{i}{2}\cdot\frac{i}{2}
\Big|
+\frac{1}{(1+s)^{1-\lambda}}
+e^{-2\int_s^t\varepsilon_ub(\tau)d\tau}
\\
\lesssim&
\frac{1}{b(s)}\Big|
\sqrt{|m_u(s,\xi)|}-\frac{1}{2}b(s)
\Big|
+\frac{1}{(1+s)^{1-\lambda}}
+e^{-2\int_s^t\varepsilon_ub(\tau)d\tau}
\\
\lesssim&
\frac{|\xi|^2+|b'(s)|}{b^2(s)}
+\frac{1}{(1+s)^{1-\lambda}}
+e^{-2\int_s^t\varepsilon_ub(\tau)d\tau}.
\end{align*}
It follows that
\begin{align*}
\Big|
\sqrt{|m_u(s,\xi)|}[\tilde{\mathscr{E}}(t,s,\xi)]_{11}
+i\frac{b(s)}{2}[\tilde{\mathscr{E}}(t,s,\xi)]_{12}
\Big|
&\lesssim
b(s)\cdot\Big|
\frac{\sqrt{|m_u(s,\xi)|}}{b(s)}[\tilde{\mathscr{E}}(t,s,\xi)]_{11}
+\frac{i}{2}[\tilde{\mathscr{E}}(t,s,\xi)]_{12}\Big|
\\
&\lesssim
\frac{|\xi|^2}{b(s)}+\frac{1}{1+s}
+\frac{1}{(1+s)^{\lambda}}e^{-2\int_s^t\varepsilon_ub(\tau)d\tau}.
\end{align*}
Compared with
$$
\Big|\sqrt{|m_u(s,\xi)|}[\tilde{\mathscr{E}}(t,s,\xi)]_{11}\Big|
\approx b(t),
\qquad
\Big|i\frac{b(s)}{2}[\tilde{\mathscr{E}}(t,s,\xi)]_{12}\Big|
\approx b(t),
$$
the multiplier $\frac{|\xi|^2}{b^2(s)}$ leads to a decay estimate multiplied by
$$
\Theta^2(t,s)\cdot(1+s)^{2\lambda},
$$
and the multiplier
$$
\frac{1}{(1+s)^{1-\lambda}}+e^{-2\int_s^t\varepsilon_ub(\tau)d\tau}
\lesssim
\frac{1}{(1+s)^{1-\lambda}}
+e^{-2\varepsilon_u((1+t)^{1-\lambda}-(1+s)^{1-\lambda})}
\lesssim \frac{1}{(1+s)^{1-\lambda}}
+(1+(1+t)^{1+\lambda}-(1+s)^{1+\lambda})^{-1},
$$
since $(1+(1+t)^{1+\lambda}-(1+s)^{1+\lambda})
e^{-2\varepsilon_u((1+t)^{1-\lambda}-(1+s)^{1-\lambda})}\lesssim 1$
for all $0\le s\le t$.
$\hfill\Box$

\begin{remark}
If $\hat C=0$, then \eqref{eq-linear-u-st-can} is reduced to
the optimal decay estimate
\eqref{eq-linear-u-opt} with the higher order regularity
$\big\|\partial_x^{|\alpha|+1+\omega_{r,q}}(v(s,\cdot),u(s,\cdot))\big\|_{L^{r}}^h$
replaced by
$\big\|\partial_x^{|\alpha|+\omega_{r,q}} u(s,\cdot)\big\|_{L^{r}}^h$.
That is, \eqref{eq-linear-u-st-can} is stronger than both \eqref{eq-linear-u}
and \eqref{eq-linear-u-opt} if $\hat C=0$.
Here we cannot prove that $\hat C=0$ due to the approximation error
in \eqref{eq-zapprox}.
Fortunately, the strategy of applying
\eqref{eq-linear-u-opt} and \eqref{eq-linear-u}
to $\|\partial_x^\alpha u(t,\cdot)\|_{L^q}$ with different index $\alpha$
works for $n\ge2$ and $\lambda\in[0,1)$.
\end{remark}

\section{Reformulated Euler system}

We apply the optimal decay estimates Theorem \ref{th-linear} of
the linear system \eqref{eq-vu} to the study of asymptotic behavior of
nonlinear system \eqref{eq-vbdu}.
We rewrite \eqref{eq-vbdu} as
\begin{equation} \label{eq-non-sys}
\partial_t
\begin{pmatrix}
v \\
\boldsymbol u
\end{pmatrix}
=
\begin{pmatrix}
0 & -\nabla\cdot \\
-\nabla & -\frac{\mu}{(1+t)^\lambda}
\end{pmatrix}
\begin{pmatrix}
v \\
\boldsymbol u
\end{pmatrix}
+
\begin{pmatrix}
-\boldsymbol u\cdot\nabla v-\varpi v\nabla\cdot \boldsymbol u \\
-(\boldsymbol u\cdot\nabla) \boldsymbol u-\varpi v\nabla v
\end{pmatrix},
\end{equation}
and the solution can be expressed as by the Duhamel principle
\begin{equation} \label{eq-Duhamel}
\begin{pmatrix}
v(t,x) \\
\boldsymbol u(t,x)
\end{pmatrix}
=\mathcal{G}(t,0)
\begin{pmatrix}
v(0,x) \\
\boldsymbol u(0,x)
\end{pmatrix}
+
\int_0^t\mathcal{G}(t,s)Q(s,x)ds,
\end{equation}
where
$$
Q(s,x)=
\begin{pmatrix}
Q_1(s,x) \\
Q_2(s,x)
\end{pmatrix}
=
\begin{pmatrix}
-\boldsymbol u\cdot\nabla v-\varpi v\nabla\cdot \boldsymbol u \\
-(\boldsymbol u\cdot\nabla) \boldsymbol u-\varpi v\nabla v
\end{pmatrix},
\quad
\mathcal{G}(t,s)
=
\begin{pmatrix}
\mathcal{G}_{11}(t,s) & \mathcal{G}_{12}(t,s) \\
\mathcal{G}_{21}(t,s) & \mathcal{G}_{22}(t,s)
\end{pmatrix}.
$$
The semigroup (Green matrix) $\mathcal{G}(t,s)$ stands for the evolution
of the linear system starting from the time $s$ to $t$.
For simplicity, we may write a function of time and space $v(t,x)$
as $v(t)$.

It should be noted that $\mathcal{G}(t,s)\ne \mathcal{G}(t-s,0)$
since the decaying damping $\frac{\mu}{(1+t)^\lambda}$ on $(s,t)$
is completely different from the damping on $(0,t-s)$.
One should be careful that the optimal decay estimates of
$\mathcal{G}(t,s)$ depends on both $t$ and $s$ (not only on $t-s$).

\subsection{Optimal $L^2$ decay estimates}

We start with the optimal $L^1$-$L^2$ decay estimates of the nonlinear system \eqref{eq-vbdu}.

\begin{lemma} \label{le-decay-G}
For $t\ge s\ge T_0$ ($T_0\ge0$ is the constant in Lemma \ref{le-Phi}), there hold
\begin{align} \nonumber
\|\partial_x^\alpha \mathcal{G}_{11}(t,s)\phi(x)\|\lesssim&
\Gamma^{\frac{n}{2}}(t,s)\cdot
\Theta^{|\alpha|}(t,s)
\cdot
(\|\phi\|_{L^1}^l+\|\partial_x^{|\alpha|}\phi\|^h),
\\ \nonumber
\|\partial_x^\alpha \mathcal{G}_{12}(t,s)\phi(x)\|\lesssim &
(1+s)^\lambda\cdot
\Gamma^{\frac{n}{2}}(t,s)\cdot
\Theta^{|\alpha|+1}(t,s)
\cdot
(\|\phi\|_{L^1}^l+\|\partial_x^{|\alpha|}\phi\|^h),
\\ \nonumber
\|\partial_x^\alpha \mathcal{G}_{21}(t,s)\phi(x)\|\lesssim &
(1+t)^\lambda\cdot
\Gamma^{\frac{n}{2}}(t,s)\cdot
\Theta^{|\alpha|+1}(t,s)
\cdot
(\|\phi\|_{L^1}^l+\|\partial_x^{|\alpha|}\phi\|^h),
\\ \label{eq-decay-G}
\|\partial_x^\alpha \mathcal{G}_{22}(t,s)\phi(x)\|\lesssim &
\Big(\frac{1+t}{1+s}\Big)^\lambda\cdot
\Gamma^{\frac{n}{2}}(t,s)\cdot
\Theta^{|\alpha|}(t,s)
\cdot
(\|\phi\|_{L^1}^l+\|\partial_x^{|\alpha|}\phi\|^h).
\end{align}

Furthermore,
\begin{align} \label{eq-decay-G-opt}
\|\partial_x^\alpha \mathcal{G}_{22}(t,s)\phi(x)\|\lesssim&
(1+t)^\lambda(1+s)^\lambda\cdot
\Gamma^{\frac{n}{2}}(t,s)\cdot
\Theta^{|\alpha|+2}(t,s)
\cdot
(\|\phi\|_{L^1}^l+\|\partial_x^{|\alpha|+1}\phi\|^h).
\end{align}
\end{lemma}
{\it\bfseries Proof.}
These estimates are simple conclusions of Theorem \ref{th-linear}.
$\hfill\Box$

\begin{lemma} \label{le-min}
For $\beta>0$ and $\gamma>0$, there holds
\begin{align} \nonumber
&\int_0^t
(1+(1+t)^{1+\lambda}-(1+s)^{1+\lambda})^{-\beta}
(1+s)^{-\gamma}ds
\\ \label{eq-min}
\lesssim&
\begin{cases}
(1+t)^{-\min\{\beta(1+\lambda),\gamma\}},
\quad&\text{if}\quad\max\{\beta(1+\lambda),\gamma\}>1,\\
(1+t)^{-\min\{\beta(1+\lambda),\gamma\}}\ln(e+t),
\quad&\text{if}\quad\max\{\beta(1+\lambda),\gamma\}=1,\\
(1+t)^{-\gamma-\beta(1+\lambda)+1},
\quad&\text{if}\quad\max\{\beta(1+\lambda),\gamma\}<1.
\end{cases}
\end{align}
\end{lemma}
{\it\bfseries Proof.}
Denote $\delta_{p,q}=1$ for $p=q$ and $\delta_{p,q}=0$ for $p\ne q$,
we can calculate
\begin{align*}
&\int_0^t
(1+(1+t)^{1+\lambda}-(1+s)^{1+\lambda})^{-\beta}
(1+s)^{-\gamma}ds
\\
\lesssim &
\Big(\int_0^{t/2}+\int_{t/2}^t\Big)
(1+(1+t)^{1+\lambda}-(1+s)^{1+\lambda})^{-\beta}
(1+s)^{-\gamma}ds
\\
\lesssim &
\int_0^{t/2}
(1+t)^{-\beta(1+\lambda)}(1+s)^{-\gamma}ds
+
\int_{t/2}^t
(1+(1+t)^{1+\lambda}-(1+s)^{1+\lambda})^{-\beta}
(1+t)^{-\gamma}ds
\\
\lesssim &
(1+t)^{-\beta(1+\lambda)}(1+t)^{\max\{1-\gamma,0\}}|\ln(e+t)|^{\delta_{\gamma,1}}
+(1+t)^{-\gamma}(1+t)^{\max\{1-\beta(1+\lambda),0\}}|\ln(e+t)|^{\delta_{\beta(1+\lambda),1}},
\end{align*}
since
\begin{align*}
&\int_{t/2}^t
(1+(1+t)^{1+\lambda}-(1+s)^{1+\lambda})^{-\beta}ds
\lesssim
\int_{t/2}^t
(1+(t-s)^{1+\lambda})^{-\beta}ds
\\
&\quad\quad\lesssim
\int_{t/2}^t
(1+t-s)^{-\beta(1+\lambda)}ds
=\int_0^{t/2}
(1+s)^{-\beta(1+\lambda)}ds.
\end{align*}
We can verify \eqref{eq-min} in different cases.
$\hfill\Box$

\begin{lemma} \label{le-min2}
For $\beta>0$, $\gamma>0$, and $k\ge0$, there holds
\begin{align} \nonumber
&\int_0^t
(1+s)^\lambda\cdot\Gamma^\beta(t,s)\cdot\Theta^{k+1}(t,s)\cdot(1+s)^{-\gamma}ds
\\ \nonumber
\lesssim&
\int_0^t
\Gamma^\beta(t,s)\cdot\Theta^{k}(t,s)\cdot(1+s)^{-\gamma}ds
\lesssim
\int_0^t
\Gamma^{\beta+k}(t,s)\cdot(1+s)^{-\gamma}ds
\\ \label{eq-min2}
\lesssim&
\begin{cases}
(1+t)^{-\min\{\frac{1+\lambda}{2}(\beta+k),\gamma\}},
\quad&\text{if}\quad\max\{\frac{1+\lambda}{2}(\beta+k),\gamma\}>1,\\
(1+t)^{-\min\{\frac{1+\lambda}{2}(\beta+k),\gamma\}}\ln(e+t),
\quad&\text{if}\quad\max\{\frac{1+\lambda}{2}(\beta+k),\gamma\}=1,\\
(1+t)^{-\gamma-\frac{1+\lambda}{2}(\beta+k)+1},
\quad&\text{if}\quad\max\{\frac{1+\lambda}{2}(\beta+k),\gamma\}<1.
\end{cases}
\end{align}
\end{lemma}
{\it\bfseries Proof.}
We note that $\Theta(t,s)=\min\{\Gamma(t,s),(1+t)^{-\lambda}\}$
as defined in \eqref{eq-Gamma}.
The proof is completed according to Lemma \ref{le-min}.
$\hfill\Box$

The following higher order energy estimates will be used to close the
decay estimates of nonlinear system \eqref{eq-vbdu}.

\begin{lemma} \label{le-energy}
Assume that $(v_0,\boldsymbol u_0)\in H^{[\frac{n}{2}]+3}$
and a priori assume that
\begin{equation} \label{eq-apriori-energy}
\|(v(t),\boldsymbol u(t))\|_{H^{[\frac{n}{2}]+2}}\le \delta_0b(t),
\end{equation}
where $\delta_0>0$ is a small constant.
Then the nonlinear system \eqref{eq-vbdu} admits a global solution $(v,\boldsymbol u)$
such that
\begin{equation} \label{eq-energy-high}
\|(v,\boldsymbol u)\|_{H^{[\frac{n}{2}]+3}}^2
+\int_0^t b(s)\big(\|\nabla v(s)\|_{H^{[\frac{n}{2}]+2}}^2
+\|\boldsymbol u(s)\|_{H^{[\frac{n}{2}]+3}}^2\big)ds
\lesssim \|(v_0,\boldsymbol u_0)\|_{H^{[\frac{n}{2}]+3}}^2.
\end{equation}
\end{lemma}
{\it\bfseries Proof.}
The energy estimate \eqref{eq-energy-high} is proved through the following four steps.
The case of time independent damping and $n=3$ is proved in \cite{TanZ-JDE13}.
Here the main difficulty lies in the absence of uniform lower bound
of the weak damping coefficient.

Step I: For $0\le k\le [\frac{n}{2}]+2$, we have
\begin{equation} \label{eq-zenergy1}
\frac{d}{dt}\|\partial_x^k(v,\boldsymbol u)\|^2
+b(t)\|\partial_x^k \boldsymbol u\|^2
\lesssim \|(v,\boldsymbol u)\|_{H^{[\frac{n}{2}]+2}}
\cdot (\|\partial_x^{k+1}v\|^2+\|\partial_x^k\boldsymbol u\|^2).
\end{equation}
This is proved by applying $\partial_x^k$ to \eqref{eq-vbdu}
and then multiplying the equation by $\partial_x^k(v,\boldsymbol u)$,
summing up and integrating over $\mathbb R^n$.
Here we omit the details.

Step II: For $0\le k\le [\frac{n}{2}]+2$, we have
\begin{equation} \label{eq-zenergy2}
\frac{d}{dt}\|\partial_x^{k+1}(v,\boldsymbol u)\|^2
+b(t)\|\partial_x^{k+1} \boldsymbol u\|^2
\lesssim \|(v,\boldsymbol u)\|_{H^{[\frac{n}{2}]+2}}
\cdot (\|\partial_x^{k+1}v\|^2+\|\partial_x^{k+1}\boldsymbol u\|^2).
\end{equation}
This is proved by applying $\partial_x^{k+1}$ to \eqref{eq-vbdu}
and then multiplying the equation by $\partial_x^{k+1}(v,\boldsymbol u)$,
summing up and integrating over $\mathbb R^n$.

Step III: For $0\le k\le [\frac{n}{2}]+2$, we have
\begin{equation} \label{eq-zenergy3}
\frac{d}{dt}\int\partial_x^k\boldsymbol u\cdot \nabla\partial_x^{k}v
+\|\partial_x^{k+1}v\|^2
\lesssim \|\partial_x^{k}\boldsymbol u\|^2
+\|(v,\boldsymbol u)\|_{H^{[\frac{n}{2}]+2}}
\cdot (\|\partial_x^{k+1}v\|^2+\|\partial_x^{k+1}\boldsymbol u\|^2).
\end{equation}
This is proved by applying $\partial_x^k$ to \eqref{eq-vbdu}$_2$
and then multiplying it by $\nabla\partial_x^{k}v$,
utilizing \eqref{eq-vbdu}$_1$ to dealing with the mixed time derivative term
$\int \partial_x^k\partial_t\boldsymbol u\cdot \nabla\partial_x^{k}v$.

Step IV: Multiply \eqref{eq-zenergy3} by $b(t)$,
for $0\le k\le [\frac{n}{2}]+2$, we have
\begin{align*}
&\frac{d}{dt}\Big(b(t)\int\partial_x^k\boldsymbol u\cdot \nabla\partial_x^{k}v\Big)
+b(t)\|\partial_x^{k+1}v\|^2
\\
\lesssim&
|b'(t)|\int\big|\partial_x^k\boldsymbol u\cdot \nabla\partial_x^{k}v\big|
+b(t)\|\partial_x^{k}\boldsymbol u\|^2
+b(t)\|(v,\boldsymbol u)\|_{H^{[\frac{n}{2}]+2}}
\cdot (\|\partial_x^{k+1}v\|^2+\|\partial_x^{k+1}\boldsymbol u\|^2)
\\
\lesssim&
\varepsilon_1b(t)\|\partial_x^{k+1}v\|^2
+b(t)\|\partial_x^{k}\boldsymbol u\|^2
+b(t)\|(v,\boldsymbol u)\|_{H^{[\frac{n}{2}]+2}}
\cdot (\|\partial_x^{k+1}v\|^2+\|\partial_x^{k+1}\boldsymbol u\|^2),
\end{align*}
where $\varepsilon_1>0$ is a small constant.
Therefore, for $0\le k\le [\frac{n}{2}]+2$,
\begin{equation} \label{eq-zenergy4}
\frac{d}{dt}\Big(b(t)\int\partial_x^k\boldsymbol u\cdot \nabla\partial_x^{k}v\Big)
+b(t)\|\partial_x^{k+1}v\|^2
\lesssim
b(t)\|\partial_x^{k}\boldsymbol u\|^2
+b(t)\|(v,\boldsymbol u)\|_{H^{[\frac{n}{2}]+2}}
\cdot (\|\partial_x^{k+1}v\|^2+\|\partial_x^{k+1}\boldsymbol u\|^2).
\end{equation}
Multiply \eqref{eq-zenergy4} by a small constant $\varepsilon_2>0$,
summing it up with \eqref{eq-zenergy1} and \eqref{eq-zenergy2},
we have
$$
\frac{d}{dt}\|(v,\boldsymbol u)\|_{H^{[\frac{n}{2}]+3}}^2
+\frac{d}{dt}\Big(\varepsilon_2\sum_{k=0}^{[{n}/{2}]+2}
b(t)\int\partial_x^k\boldsymbol u\cdot\nabla \partial_x^{k}v\Big)
+b(t)(\|\nabla v\|_{H^{[\frac{n}{2}]+2}}^2
+\|\boldsymbol u\|_{H^{[\frac{n}{2}]+3}}^2)\le0,
$$
provided that the a priori assumption \eqref{eq-apriori-energy} is valid.
The constant $\varepsilon_2>0$ is small such that
$$
\Big|\varepsilon_2\sum_{k=0}^{[{n}/{2}]+2}
b(t)\int\partial_x^k\boldsymbol u\cdot \nabla\partial_x^{k}v\Big|
\le \frac{1}{2}\|(v,\boldsymbol u)\|_{H^{[\frac{n}{2}]+3}}^2.
$$
The proof is completed.
$\hfill\Box$

We present the optimal $L^1$-$L^2$ decay rates of the nonlinear system \eqref{eq-vbdu}.

\begin{proposition}[Decay rates of nonlinear system] \label{th-nonlinear-n4}
For $n\ge2$ and $\lambda\in[0,1)$, there exists a constant $\varepsilon_0>0$, such that
the solution $(v,\boldsymbol u)$ of the nonlinear system \eqref{eq-vbdu}
corresponding to initial data $(v_0,\boldsymbol u_0)$
with small energy $\|(v_0,\boldsymbol u_0)\|_{L^1\cap H^{[\frac{n}{2}]+3}}\le\varepsilon_0$
exists globally and satisfies
\begin{equation} \label{eq-nonlinear-n4}
\begin{cases}
\|\partial_x^\alpha v\|\lesssim (1+t)^{-\frac{1+\lambda}{4}n-\frac{1+\lambda}{2}|\alpha|},
\quad &0\le |\alpha|\le [\frac{n}{2}]+1,\\
\|\partial_x^\alpha \boldsymbol u\|\lesssim
(1+t)^{-\frac{1+\lambda}{4}n-\frac{1+\lambda}{2}(|\alpha|+1)+\lambda},
\quad &0\le |\alpha|\le [\frac{n}{2}],\\
\|\partial_x^\alpha \boldsymbol u\|\lesssim
(1+t)^{-\frac{1+\lambda}{4}n-\frac{1+\lambda}{2}|\alpha|+\lambda},
\quad &|\alpha|=[\frac{n}{2}]+1, \\
\|(v,\boldsymbol u)\|_{H^{[\frac{n}{2}]+3}}\lesssim 1.
\end{cases}
\end{equation}
The first two decay estimates in \eqref{eq-nonlinear-n4}
(i.e., the decay estimates on $\|\partial_x^\alpha v\|$
with $0\le|\alpha|\le [\frac{n}{2}]+1$ and
$\|\partial_x^\alpha \boldsymbol u\|$ with $0\le|\alpha|\le [\frac{n}{2}]$) are optimal.
\end{proposition}
{\it\bfseries Proof.}
Suppose that the local solution $(v,\boldsymbol u)$ exists for $t\in(0,T)$.
Denote the weighted energy
\begin{align} \nonumber
E_n(\tilde t)&:=\sup_{t\in(0,\tilde t)}\Big\{
\sum_{0\le|\alpha|\le [n/2]+1}
(1+t)^{\frac{1+\lambda}{4}n+\frac{1+\lambda}{2}|\alpha|}\|\partial_x^\alpha v\|,
\sum_{0\le|\alpha|\le [n/2]}
(1+t)^{\frac{1+\lambda}{4}n+\frac{1+\lambda}{2}(|\alpha|+1)-\lambda}
\|\partial_x^\alpha \boldsymbol u\|,
\\
&\sum_{|\alpha|=[n/2]+1}
(1+t)^{\frac{1+\lambda}{4}n+\frac{1+\lambda}{2}|\alpha|-\lambda}
\|\partial_x^\alpha \boldsymbol u\|,
\sum_{|\alpha|= [n/2]+2}
(1+t)^{\frac{1+\lambda}{4}n}
\|\partial_x^\alpha (v,\boldsymbol u)\|,
\sum_{|\alpha|=[n/2]+3}
\|\partial_x^\alpha (v,\boldsymbol u)\|
\Big\}.
\end{align}
We claim that under the condition
$\|(v_0,\boldsymbol u_0)\|_{L^1\cap H^{[\frac{n}{2}]+3}}\le\varepsilon_0$,
there holds
\begin{equation} \label{eq-apriori-n4}
E_n(\tilde t)\lesssim \delta_0, \quad \forall \tilde t\in(0,T),
\end{equation}
where $\varepsilon_0>0$ and $\delta_0>0$ are small constants to be determined.

The global existence and the a priori assumption \eqref{eq-apriori-n4}
(which implies the decay estimates \eqref{eq-nonlinear-n4})
will be proved in the following three steps.
For the sake of simplicity, we take the case $n=3$ for example.
Other cases with $n\ge2$ follow similarly.
We may assume that $T_0=0$, where $T_0\ge0$ is the constant in Lemma \ref{le-Phi}.
That is, we consider the nonlinear system \eqref{eq-vbdu} starting form the time
$T_0$ and we write $t-T_0$ as $t$ for convenience.

Step I: Basic energy decay estimates.

According to the Duhamel principle \eqref{eq-Duhamel}
and the decay estimates of the Green matrix $\mathcal{G}(t,s)$
in Lemma \ref{le-decay-G}, we have
\begin{align*}
\|v(t)\|
\lesssim&
\|\mathcal{G}_{11}(t,0)v_0\|+\|\mathcal{G}_{12}(t,0)\boldsymbol u_0\|
+\int_0^t\|\mathcal{G}_{11}(t,s)Q_1(s)\|ds
+\int_0^t\|\mathcal{G}_{12}(t,s)Q_2(s)\|ds
\\
\lesssim&
\varepsilon_0(1+t)^{-\frac{1+\lambda}{4}n}
+\int_0^t
\Gamma^{\frac{n}{2}}(t,s)\cdot
(\|Q_1(s)\|_{L^1}^l+\|Q_1(s)\|^h)ds
\\
&+\int_0^t
(1+s)^\lambda\cdot
\Gamma^{\frac{n}{2}}(t,s)
\cdot
\Theta(t,s)
\cdot
(\|Q_2(s)\|_{L^1}^l+\|Q_2(s)\|^h)ds
\\
\lesssim&
\varepsilon_0(1+t)^{-\frac{1+\lambda}{4}n}
+E_n^2(t)\int_0^t
\Gamma^{\frac{n}{2}}(t,s)\cdot
(1+s)^{-\frac{1+\lambda}{2}n-1}ds
\\
&+E_n^2(t)\int_0^t
(1+s)^\lambda\cdot
\Gamma^{\frac{n}{2}}(t,s)
\cdot
\Theta(t,s)
\cdot
(1+s)^{-\frac{1+\lambda}{2}n-\frac{1+\lambda}{2}}ds
\\
\lesssim&
\varepsilon_0(1+t)^{-\frac{1+\lambda}{4}n}
+E_n^2(t)(1+t)^{-\frac{1+\lambda}{4}n},
\end{align*}
where we have used Lemma \ref{le-min2}
(note that $\frac{1+\lambda}{2}n+\frac{1+\lambda}{2}>1$
for $n\ge2$ and $\lambda\in[0,1)$) and
the following decay estimates on $\|Q(s)\|_{L^1}$ and $\|Q(s)\|$
(here and after, we use $D^j:=\partial_x^j$
and we may also write $\boldsymbol u$ as $u$ for simplicity)
\begin{align*}
\|Q_1(s)\|_{L^1}
&\lesssim \|uDv\|_{L^1}+\|vDu\|_{L^1}
\lesssim \|u\|\|Dv\|+\|v\|\|Du\|
\lesssim E_n^2(s)(1+s)^{-\frac{1+\lambda}{2}n-1},\\
\|Q_2(s)\|_{L^1}
&\lesssim \|uDu\|_{L^1}+\|vDv\|_{L^1}
\lesssim \|u\|\|Du\|+\|v\|\|Dv\|
\lesssim E_n^2(s)(1+s)^{-\frac{1+\lambda}{2}n-\frac{1+\lambda}{2}}.
\end{align*}
For $n=3$, we have
\begin{align*}
\|u(s)\|_{L^\infty}
&\lesssim \|Du\|^\frac{1}{2}\|D^2u\|^\frac{1}{2}
\lesssim E_n(s)(1+s)^{-\frac{1+\lambda}{4}n-1},\\
\|v(s)\|_{L^\infty}
&\lesssim \|Dv\|^\frac{1}{2}\|D^2v\|^\frac{1}{2}
\lesssim E_n(s)(1+s)^{-\frac{1+\lambda}{4}n-\frac{1+\lambda}{2}\cdot3},\\
\|Du(s)\|_{L^\infty}
&\lesssim \|D^2u\|^\frac{1}{2}\|D^3u\|^\frac{1}{2}
\lesssim E_n(s)(1+s)^{-\frac{1+\lambda}{4}n-\frac{1}{2}}, \\
\|Dv(s)\|_{L^\infty}
&\lesssim \|D^2v\|^\frac{1}{2}\|D^3v\|^\frac{1}{2}
\lesssim E_n(s)(1+s)^{-\frac{1+\lambda}{4}n-\frac{1}{2}(1+\lambda)},\\
\|D^2u(s)\|_{L^\infty}
&\lesssim \|D^3u\|^\frac{1}{2}\|D^4u\|^\frac{1}{2}
\lesssim E_n(s)(1+s)^{-\frac{1+\lambda}{8}n}, \\
\|D^2v(s)\|_{L^\infty}
&\lesssim \|D^3v\|^\frac{1}{2}\|D^4v\|^\frac{1}{2}
\lesssim E_n(s)(1+s)^{-\frac{1+\lambda}{8}n},
\end{align*}
and
\begin{align*}
\|Q_1(s)\|
&\lesssim \|uDv\|+\|vDu\|
\lesssim \|u\|_{L^\infty}\|Dv\|+\|v\|_{L^\infty}\|Du\|
\lesssim E_n^2(s)(1+s)^{-\frac{1+\lambda}{2}n-\frac{1+\lambda}{2}-1}, \\
\|Q_2(s)\|
&\lesssim \|uDu\|+\|vDv\|
\lesssim \|u\|_{L^\infty}\|Du\|+\|v\|_{L^\infty}\|Dv\|
\lesssim E_n^2(s)(1+s)^{-\frac{1+\lambda}{2}n-2},\\
\|DQ_1(s)\|
&\lesssim \|DuDv\|+\|uD^2v\|+\|vD^2u\|
\lesssim E_n^2(s)(1+s)^{-\frac{1+\lambda}{2}n-1-\frac{1+\lambda}{2}}, \\
\|DQ_2(s)\|
&\lesssim \|uD^2u\|+\|DuDu\|+\|vD^2v\|+\|DvDv\|
\lesssim E^2(s)(1+s)^{-\frac{1+\lambda}{2}n-\theta_{12}},
\end{align*}
where $\theta_{12}=\min\{\frac{3}{2},1+\lambda\}\ge \frac{1+\lambda}{2}$.

Using the above estimates, we have
\begin{align*}
\|Dv(t)\|
\lesssim&
\|D\mathcal{G}_{11}(t,0)v_0\|+\|D\mathcal{G}_{12}(t,0)u_0\|
+\int_0^t\|D\mathcal{G}_{11}(t,s)Q_1(s)\|ds
+\int_0^t\|D\mathcal{G}_{12}(t,s)Q_2(s)\|ds
\\
\lesssim&
\varepsilon_0(1+t)^{-\frac{1+\lambda}{4}n-\frac{1+\lambda}{2}}
+\int_0^t
\Gamma^{\frac{n}{2}}(t,s)
\cdot
\Theta(t,s)
\cdot
(\|Q_1(s)\|_{L^1}+\|DQ_1(s)\|)ds
\\
&+\int_0^t
(1+s)^\lambda\cdot
\Gamma^{\frac{n}{2}}(t,s)
\cdot
\Theta^2(t,s)
\cdot
(\|Q_2(s)\|_{L^1}+\|DQ_2(s)\|)ds
\\
\lesssim&
\varepsilon_0(1+t)^{-\frac{1+\lambda}{4}n-\frac{1+\lambda}{2}}
+E_n^2(t)
\int_0^t
\Gamma^{\frac{n}{2}}(t,s)
\cdot
\Theta(t,s)
\cdot
(1+s)^{-\frac{1+\lambda}{2}n-1}ds
\\
&+E_n^2(t)\int_0^t
(1+s)^\lambda\cdot
\Gamma^{\frac{n}{2}}(t,s)
\cdot
\Theta^2(t,s)
\cdot
(1+s)^{-\frac{1+\lambda}{2}n-\frac{1+\lambda}{2}}ds
\\
\lesssim&
\varepsilon_0(1+t)^{-\frac{1+\lambda}{4}n-\frac{1+\lambda}{2}}
+E_n^2(t)(1+t)^{-\frac{1+\lambda}{4}n-\frac{1+\lambda}{2}},
\end{align*}
and
\begin{align*}
\|D^2v(t)\|
\lesssim&
\|D^2\mathcal{G}_{11}(t,0)v_0\|+\|D^2\mathcal{G}_{12}(t,0)u_0\|
+\int_0^t\|D^2\mathcal{G}_{11}(t,s)Q_1(s)\|ds
+\int_0^t\|D^2\mathcal{G}_{12}(t,s)Q_2(s)\|ds
\\
\lesssim&
\varepsilon_0(1+t)^{-\frac{1+\lambda}{4}n-(1+\lambda)}
+\int_0^t
\Gamma^{\frac{n}{2}}(t,s)
\cdot
\Theta^2(t,s)
\cdot
(\|Q_1(s)\|_{L^1}+\|D^2Q_1(s)\|)ds
\\
&+\int_0^t
(1+s)^\lambda\cdot
\Gamma^{\frac{n}{2}}(t,s)
\cdot
\Theta^3(t,s)
\cdot
(\|Q_2(s)\|_{L^1}+\|D^2Q_2(s)\|)ds
\\
\lesssim&
\varepsilon_0(1+t)^{-\frac{1+\lambda}{4}n-(1+\lambda)}
+E_n^2(t)\int_0^t
\Gamma^{\frac{n}{2}}(t,s)
\cdot
\Theta^2(t,s)
\cdot
(1+s)^{-\frac{1+\lambda}{2}n-1}ds
\\
&+E_n^2(t)\int_0^t
(1+s)^\lambda\cdot
\Gamma^{\frac{n}{2}}(t,s)
\cdot
\Theta^3(t,s)
\cdot
(1+s)^{-\frac{1+\lambda}{2}n-\frac{1+\lambda}{2}}ds
\\
\lesssim&
\varepsilon_0(1+t)^{-\frac{1+\lambda}{4}n-(1+\lambda)}
+E_n^2(t)(1+t)^{-\frac{1+\lambda}{4}n-(1+\lambda)},
\end{align*}
since
$$
\frac{1+\lambda}{2}n+\frac{1+\lambda}{2}\ge \frac{1+\lambda}{4}n+(1+\lambda),
\quad \text{~for ~} n\ge2, \lambda\in[0,1).
$$
We have also used the following estimates
\begin{align*}
\|D^2Q_1(s)\|
&\lesssim \|uD^3v\|+\|DuD^2v\|+\|DvD^2u\|+\|vD^3u\|
\lesssim E_n^2(s)(1+s)^{-\frac{1+\lambda}{2}n-1}, \\
\|D^2Q_2(s)\|
&\lesssim \|uD^3u\|+\|DuD^2u\|+\|vD^3v\|+\|DvD^2v\|
\lesssim E_n^2(s)(1+s)^{-\frac{1+\lambda}{2}n-1}.
\end{align*}

The decay estimates on
$\|\partial_x^\alpha v\|$ for $0\le|\alpha|\le [\frac{n}{2}]+1$
are based on the optimal decay estimates on
$\|\partial_x^\alpha G_{11}(t,s)\|$ and $\|\partial_x^\alpha G_{12}(t,s)\|$
in \eqref{eq-decay-G}.
However, the estimates on
$\|\partial_x^\alpha G_{21}(t,s)\|$ and $\|\partial_x^\alpha G_{22}(t,s)\|$
in \eqref{eq-decay-G} is insufficient for the optimal decay estimates
on $\|\partial_x^\alpha \boldsymbol u\|$ for $0\le|\alpha|\le [\frac{n}{2}]$.
In fact, we use the optimal decay estimates in \eqref{eq-decay-G-opt}
to show the decay estimates on
$\|\partial_x^\alpha \boldsymbol u\|$ for $0\le|\alpha|\le [\frac{n}{2}]$
in a similar way as $\|\partial_x^\alpha v\|$ for $1\le|\alpha|\le [\frac{n}{2}]+1$.
One can check that the condition on the estimate of
$\|\partial_x^k \boldsymbol u\|$ for $0\le k\le [\frac{n}{2}]$
is equivalent to the condition on the estimate of
$\|\partial_x^{k+1} v\|$.

Further, we use the decay estimates in \eqref{eq-decay-G}
to show the decay estimates on
$\|\partial_x^\alpha \boldsymbol u\|$ for $[\frac{n}{2}]+1\le|\alpha|\le [\frac{n}{2}]+2$
since the regularity required in \eqref{eq-decay-G}
is one order lower than that in \eqref{eq-decay-G-opt}.
We note that in this case the condition on the estimate of
$\|\partial_x^k \boldsymbol u\|$ for $[\frac{n}{2}]+1\le k\le [\frac{n}{2}]+2$
is similar to the condition on the estimate of $\|\partial_x^k v\|$.
We have
\begin{align*}
\|D^3Q(s)\|
&\lesssim \|(v,u)D^4(v,u)\|+\|D(v,u)D^3(v,u)\|+\|D^2(v,u)D^2(v,u)\|
\lesssim E_n^2(s)(1+s)^{-\frac{1+\lambda}{4}n-\theta_3},
\end{align*}
with
\begin{align*}
\theta_3&=\min\Big\{
\frac{1+\lambda}{2}\cdot3,
\frac{1+\lambda}{4}n+\frac{1}{2},
1+\frac{1+\lambda}{8}n\Big\}.
\end{align*}
Therefore,
\begin{align*}
\|D^3u(t)\|
\lesssim&
\|D^3\mathcal{G}_{21}(t,0)v_0\|+\|D^3\mathcal{G}_{22}(t,0)u_0\|
+\int_0^t\|D^3\mathcal{G}_{21}(t,s)Q_1(s)\|ds
+\int_0^t\|D^3\mathcal{G}_{22}(t,s)Q_2(s)\|ds
\\
\lesssim&
\varepsilon_0(1+t)^{-\frac{1+\lambda}{4}n-\frac{3}{2}(1+\lambda)+\lambda}
\\
&+\int_0^t(1+t)^\lambda
\cdot
\Gamma^{\frac{n}{2}}(t,s)
\cdot
\Theta^4(t,s)
\cdot
(\|Q_1(s)\|_{L^1}+\|D^3Q_1(s)\|)ds
\\
&+\int_0^t
\Big(\frac{1+t}{1+s}\Big)^\lambda
\cdot
\Gamma^{\frac{n}{2}}(t,s)
\cdot
\Theta^3(t,s)
\cdot
(\|Q_2(s)\|_{L^1}+\|D^3Q_2(s)\|)ds
\\
\lesssim&
\varepsilon_0(1+t)^{-\frac{1+\lambda}{4}n-\frac{3}{2}(1+\lambda)+\lambda}
+E_n^2(t)\int_0^t(1+t)^\lambda
\cdot
\Gamma^{\frac{n}{2}}(t,s)
\cdot
\Theta^4(t,s)
\cdot
(1+s)^{-\min\{\frac{1+\lambda}{4}n+\theta_3,\frac{1+\lambda}{2}n+1\}}ds
\\
&+E_n^2(t)\int_0^t
\Big(\frac{1+t}{1+s}\Big)^\lambda
\cdot
\Gamma^{\frac{n}{2}}(t,s)
\cdot
\Theta^3(t,s)
\cdot
(1+s)^{-\min\{\frac{1+\lambda}{4}n+\theta_3,\frac{1+\lambda}{2}n+\frac{1+\lambda}{2}\}}ds
\\
\lesssim&
\varepsilon_0(1+t)^{-\frac{1+\lambda}{4}n}
+E_n^2(t)(1+t)^{-\frac{1+\lambda}{4}n},
\end{align*}
since $\theta_3\ge\lambda$.
The estimates on $\|D^3v\|$ follows similarly.

Step II: Higher order energy estimates.
We note that the condition \eqref{eq-apriori-n4} is stronger than
the a priori assumption \eqref{eq-apriori-energy},
and according to \eqref{eq-energy-high} in Lemma \ref{le-energy},
we have
\begin{equation} \label{eq-zenergy-high}
\|(v,\boldsymbol u)\|_{H^{[\frac{n}{2}]+3}}^2
+\int_0^t b(s)\big(\|\nabla v(s)\|_{H^{[\frac{n}{2}]+2}}^2
+\|\boldsymbol u(s)\|_{H^{[\frac{n}{2}]+3}}^2\big)ds
\lesssim \|(v_0,\boldsymbol u_0)\|_{H^{[\frac{n}{2}]+3}}^2.
\end{equation}

Step III: Closure of the a priori estimate \eqref{eq-apriori-n4}.
Combining the above estimates and choosing $\varepsilon_0>0$ and $\delta_0>0$
to be sufficiently small such that
$C(\varepsilon_0+\delta_0^2)\le \delta_0$, we see that the a priori estimate
\eqref{eq-apriori-n4} holds for all the time $t\in(0,+\infty)$.

Finally, we show that those estimates
($\|\partial_x^\alpha v\|$ with $0\le|\alpha|\le[\frac{n}{2}]+1$ and
$\|\partial_x^\alpha \boldsymbol u\|$ with $0\le|\alpha|\le[\frac{n}{2}]$)
are optimal.
We take the estimate on $\|v\|$ for example.
According to the optimal decay estimates Lemma \ref{le-decay-G} and
the energy estimates in Step I,
we choose the initial data $(v_0,\boldsymbol u_0)$ such that
$\|\mathcal{G}_{11}(t,0)v_0\|$ decays optimally, then we have
$$
\|v(t)\|
\gtrsim
\|\mathcal{G}_{11}(t,0)v_0\|-\|\mathcal{G}_{12}(t,0)\boldsymbol u_0\|
-\int_0^t\|\mathcal{G}_{11}(t,s)Q_1(s)\|ds
-\int_0^t\|\mathcal{G}_{12}(t,s)Q_2(s)\|ds,
$$
where $\|\mathcal{G}_{12}(t,0)\boldsymbol u_0\|$ decays faster than
$\|\mathcal{G}_{11}(t,0)v_0\|$,
and $\int_0^t\|\mathcal{G}_{11}(t,s)Q_1(s)\|ds
+\int_0^t\|\mathcal{G}_{12}(t,s)Q_2(s)\|ds$ decays
no slower than $\|\mathcal{G}_{11}(t,0)v_0\|$.
We note that $Q_1(t,x)$ and $Q_2(t,x)$ are quadratic,
and we rescale the initial data as $(\varepsilon_1v_0,\varepsilon_1\boldsymbol u_0)$
with $\varepsilon_1>0$ sufficiently small
such that neither
$\int_0^t\|\mathcal{G}_{11}(t,s)Q_1(s)\|ds$ nor
$\int_0^t\|\mathcal{G}_{12}(t,s)Q_2(s)\|ds$
is comparable with $\|\mathcal{G}_{11}(t,0)v_0\|$.
In fact, according to the proof in Step I, we have
$$
\int_0^t\|\mathcal{G}_{1j}(t,s)Q_j(s)\|ds\lesssim E_n^2(t)(1+t)^{-\frac{1+\lambda}{4}n}
\lesssim \delta_0^2(1+t)^{-\frac{1+\lambda}{4}n},
\quad j=1,2,
$$
and
$$
\|\mathcal{G}_{11}(t,0)v_0\|\approx \varepsilon_0(1+t)^{-\frac{1+\lambda}{4}n},
$$
where the small constants $\varepsilon_0\approx \delta_0$
as in the proof Step III.
That is, $\|v(t)\|$ decays in the same order as
$\|\mathcal{G}_{11}(t,0)v_0\|$.
The proof is completed.
$\hfill\Box$

{\it\bfseries Proof of Theorem \ref{th-nonlinear}.}
This is proved in Proposition \ref{th-nonlinear-n4}.
$\hfill\Box$

\subsection{Optimal $L^q$ decay estimates}

We now turn to the $L^1$-$L^q$ decay estimates of the nonlinear system \eqref{eq-vbdu}.
Similar to Lemma \ref{le-energy}, we have the following higher order energy estimates.

\begin{lemma} \label{le-energy-h}
Assume that $(v_0,\boldsymbol u_0)\in H^{[\frac{n}{2}]+k}$ with $k\ge2$
and a priori assume that
\begin{equation*}
\|(v(t),\boldsymbol u(t))\|_{H^{[\frac{n}{2}]+2}}\le \delta_0b(t),
\end{equation*}
where $\delta_0>0$ is a small constant.
Then the nonlinear system \eqref{eq-vbdu} admits a global solution $(v,\boldsymbol u)$
such that
\begin{equation} \label{eq-energy-high-h}
\|(v,\boldsymbol u)\|_{H^{[\frac{n}{2}]+k}}^2
+\int_0^t b(s)\big(\|\nabla v(s)\|_{H^{[\frac{n}{2}]+k-1}}^2
+\|\boldsymbol u(s)\|_{H^{[\frac{n}{2}]+k}}^2\big)ds
\lesssim \|(v_0,\boldsymbol u_0)\|_{H^{[\frac{n}{2}]+k}}^2.
\end{equation}
\end{lemma}
{\it\bfseries Proof.}
The proof is completely same as that in Lemma \ref{le-energy}.
We note that the a priori assumption only requires the norms
$\|(v(t),\boldsymbol u(t))\|_{H^{[\frac{n}{2}]+2}}$,
which is sufficient for the required estimates such as
$\|\partial_x(v(t),\boldsymbol u(t))\|_{L^\infty}$ and
$\|(v(t),\boldsymbol u(t))\|_{L^\infty}$.
$\hfill\Box$

\begin{lemma} \label{le-decay-G-h}
For $q\in[2,\infty]$ and $1\le p,r\le 2$
(or $\theta\in[0,\frac{n}{2})$),
and for $t\ge s\ge T_0$ ($T_0\ge0$ is the constant in Lemma \ref{le-Phi}), we have
\begin{align*}
\|\partial_x^\alpha \mathcal{G}_{11}(t,s)\phi(x)\|_{L^q}\lesssim&
\Gamma^{\gamma_{p,q}}(t,s)\cdot
\Theta^{|\alpha|}(t,s)
\cdot
(\|\phi\|_{L^{p}}^l+\|\partial_x^{|\alpha|+\omega_{r,q}}\phi\|_{L^{r}}^h),
\\
\|\partial_x^\alpha \mathcal{G}_{12}(t,s)\phi(x)\|_{L^q}\lesssim &
(1+s)^\lambda\cdot
\Gamma^{\gamma_{p,q}}(t,s)\cdot
\Theta^{|\alpha|+1}(t,s)
\cdot
(\|\phi\|_{L^{p}}^l+\|\partial_x^{|\alpha|+\omega_{r,q}}\phi\|_{L^{r}}^h),
\\
\|\partial_x^\alpha \mathcal{G}_{21}(t,s)\phi(x)\|_{L^q}\lesssim &
(1+t)^\lambda\cdot
\Gamma^{\gamma_{p,q}}(t,s)\cdot
\Theta^{|\alpha|+1}(t,s)
\cdot
(\|\phi\|_{L^{p}}^l+\|\partial_x^{|\alpha|+\omega_{r,q}}\phi\|_{L^{r}}^h),
\\
\|\partial_x^\alpha \mathcal{G}_{22}(t,s)\phi(x)\|_{L^q}\lesssim &
\Big(\frac{1+t}{1+s}\Big)^\lambda\cdot
\Gamma^{\gamma_{p,q}}(t,s)\cdot
\Theta^{|\alpha|}(t,s)
\cdot
(\|\phi\|_{L^{p}}^l+\|\partial_x^{|\alpha|+\omega_{r,q}}\phi\|_{L^{r}}^h),
\end{align*}
where $\gamma_{p,q}:=n(1/{p}-1/q)$
(or $\gamma_{p,q}$ replaced by $\beta_{\theta,q}:=\theta+\gamma_{2,q}$
and $\|\cdot\|_{L^{p}}$ norm replaced by $\|\cdot\|_{\dot H^{-\theta}}$),
and $\omega_{r,q}>\gamma_{r,q}$ for $(r,q)\ne(2,2)$ and $\omega_{2,2}=0$.

Furthermore,
\begin{align*}
\|\partial_x^\alpha \mathcal{G}_{22}(t,s)\phi(x)\|_{L^q}\lesssim&
(1+t)^\lambda(1+s)^\lambda\cdot
\Gamma^{\gamma_{p,q}}(t,s)\cdot
\Theta^{|\alpha|+2}(t,s)
\cdot
(\|\phi\|_{L^{p}}^l+\|\partial_x^{|\alpha|+1+\omega_{r,q}}\phi\|_{L^{r}}^h).
\end{align*}
\end{lemma}
{\it\bfseries Proof.}
These estimates are simple conclusions of Theorem \ref{th-linear}.
$\hfill\Box$

We present the following optimal $L^q$ decay estimates
of the nonlinear system \eqref{eq-vbdu}.

\begin{proposition}[Optimal $L^q$ decay estimates] \label{th-nonlinear-Lp-2}
For $n\ge2$, $\lambda\in[0,1)$, $q\in[2,\infty]$ and $k\ge 3+[\gamma_{2,q}]$
with $\gamma_{2,q}:=n(1/2-1/q)$,
let $(v,\boldsymbol u)$ be the solution to the nonlinear system \eqref{eq-vbdu}
corresponding to initial data $(v_0,\boldsymbol u_0)$
with small energy such that
$\|(v_0,\boldsymbol u_0)\|_{L^1\cap H^{[\frac{n}{2}]+k}}\le\varepsilon_0$,
where $\varepsilon_0>0$ is a small constant only depending on $n,q,k$ and
the constants $\gamma,\mu,\lambda$ in the system.
Then $(v,\boldsymbol u)\in L^\infty(0,+\infty;H^{[\frac{n}{2}]+k})$ and satisfies
\begin{equation} \label{eq-nonlinear-Lp-3}
\begin{cases}
\|\partial_x^\alpha v\|_{L^q}\lesssim
(1+t)^{-\frac{1+\lambda}{2}\gamma_{1,q}-\frac{1+\lambda}{2}|\alpha|},
\quad &0\le |\alpha|\le 1,\\
\|\boldsymbol u\|_{L^q}\lesssim
(1+t)^{-\frac{1+\lambda}{2}\gamma_{1,q}-\frac{1-\lambda}{2}},
\end{cases}
\end{equation}
where $\gamma_{1,q}=n(1-1/q)$.
All the decay estimates in \eqref{eq-nonlinear-Lp-3} are optimal.
\end{proposition}
{\it\bfseries Proof.}
The decay estimates
are based on the optimal $L^2$ decay estimates Proposition \ref{th-nonlinear-n4},
the higher order energy estimates Lemma \ref{le-energy-h},
and the $L^1$-$L^q$ decay estimates of the Green matrix in Lemma \ref{le-decay-G-h}.

We prove the estimate on $\|\partial_x^\alpha v\|_{L^q}$ with $|\alpha|=1$
in \eqref{eq-nonlinear-Lp-3}.
According to the Duhamel principle \eqref{eq-Duhamel} and
the $L^1$-$L^q$ decay estimates of the Green matrix in Lemma \ref{le-decay-G-h},
we have
\begin{eqnarray*}
\|Dv(t)\|_{L^q}
&\lesssim&
\|D\mathcal{G}_{11}(t,0)v_0\|_{L^q}+\|D\mathcal{G}_{12}(t,0)u_0\|_{L^q} \\
&&+\int_0^t\|D\mathcal{G}_{11}(t,s)Q_1(s)\|_{L^q}ds
+\int_0^t\|D\mathcal{G}_{12}(t,s)Q_2(s)\|_{L^q}ds
\\
&\lesssim&
\varepsilon_0(1+t)^{-\frac{1+\lambda}{2}\gamma_{1,q}-\frac{1+\lambda}{2}}
+\int_0^t
\Gamma^{\gamma_{1,q}}(t,s)\cdot
\Theta(t,s)
\cdot
(\|Q_1(s)\|_{L^1}+\|D^{1+\omega_{2,q}}Q_1(s)\|)ds
\\
&&+\int_0^t
(1+s)^\lambda\cdot
\Gamma^{\gamma_{1,q}}(t,s)\cdot
\Theta^2(t,s)
\cdot
(\|Q_2(s)\|_{L^1}+\|D^{1+\omega_{2,q}}Q_2(s)\|)ds
\\
&\lesssim&
\varepsilon_0(1+t)^{-\frac{1+\lambda}{2}\gamma_{1,q}-\frac{1+\lambda}{2}}
+E_n^2(t)\int_0^t
\Gamma^{\gamma_{1,q}}(t,s)\cdot
\Theta(t,s)
\cdot
(1+s)^{-\frac{1+\lambda}{2}n-1}ds
\\
&&+E_n^2(t)\int_0^t
(1+s)^\lambda\cdot
\Gamma^{\gamma_{1,q}}(t,s)\cdot
\Theta^2(t,s)
\cdot
(1+s)^{-\frac{1+\lambda}{2}n-\frac{1+\lambda}{2}}ds
\\
&\lesssim&
\varepsilon_0(1+t)^{-\frac{1+\lambda}{2}\gamma_{1,q}-\frac{1+\lambda}{2}}
+E_n^2(t)(1+t)^{-\frac{1+\lambda}{2}\gamma_{1,q}-\frac{1+\lambda}{2}},
\end{eqnarray*}
where $\omega_{2,q}>\gamma_{2,q}$ and
$$
\frac{1+\lambda}{2}n+\frac{1+\lambda}{2}
\ge
\frac{1+\lambda}{2}\gamma_{1,q}+\frac{1+\lambda}{2},
$$
which is valid for all $n\ge2$, $\lambda\in[0,1)$, and $q\in[2,\infty]$.
Other estimates and cases can be proved through a similar procedure.
$\hfill\Box$

{\it\bfseries Proof of Theorem \ref{th-nonlinear-Lp}.}
This is proved in Proposition \ref{th-nonlinear-Lp-2}.
$\hfill\Box$

\

{\bf Acknowledgement}.
This work was done when the first author visited McGill University
supported by China Scholarship Council (CSC) for
the senior visiting scholar program.
He would like to express his sincere thanks for the hospitality
of McGill University and CSC.
The research of the first author was supported by
NSFC Grant No.~11701184 and CSC No. 201906155021.
The research of the second was supported in part
by NSERC Grant RGPIN 354724-16, and FRQNT Grant No. 2019-CO-256440.

\end{document}